\documentclass[10pt]{article}

\usepackage[english]{babel}
\usepackage[letterpaper,top=2cm,bottom=2cm,left=3cm,right=3cm]{geometry}
\usepackage[utf8]{inputenc}
\usepackage{amsfonts,amsmath,amssymb,amsthm,mathtools}
\usepackage{algorithm}
\usepackage{algorithmic}
\usepackage{array}
\usepackage{booktabs}
\usepackage{cite}
\usepackage[dvipsnames]{xcolor}
\usepackage{graphicx}
\usepackage{float}
\usepackage[
  colorlinks=true,
  allcolors=NavyBlue,
  urlcolor=Black
]{hyperref}

\hypersetup{
  pdftitle={Accelerated Reflected Splitting with Degenerate Preconditioning for Monotone Inclusions},
  pdfauthor={Jiahong Li, Defeng Sun, Hailing Wang, Liping Zhang, Wei Zhao}
}

\theoremstyle{plain}
\newtheorem{theorem}{Theorem}[section]
\newtheorem{lemma}[theorem]{Lemma}

\theoremstyle{definition}
\newtheorem{assumption}[theorem]{Assumption}

\theoremstyle{remark}
\newtheorem{remark}[theorem]{Remark}

\newcommand{\Hcal}{\mathcal{H}}

\newcommand{\Vcal}{\mathcal V}
\newcommand{\Wcal}{\mathcal W}
\newcommand{\Scal}{\mathcal S}
\newcommand{\Id}{\operatorname{Id}}
\newcommand{\Fix}{\operatorname{Fix}}
\newcommand{\ran}{\operatorname{ran}}
\newcommand{\dist}{\operatorname{dist}}

\newcommand{\rnat}{r_{\mathrm{nat}}}
\newcommand{\norm}[1]{\left\lVert #1\right\rVert}
\newcommand{\ip}[2]{\left\langle #1,#2\right\rangle}
\newcommand{\R}{\mathbb{R}}

\title{Accelerated Reflected Splitting with Degenerate Preconditioning for Monotone Inclusions\thanks{The research of Defeng Sun was supported in part by the Hong Kong RGC Senior Research Fellow Scheme under Grant No. SRFS2223-5S02 and by the Research Centre for Intelligent Operations Research.
The research of Liping Zhang was supported in part by the National Natural Science Foundation of China under Grant No. 12571323.}}
\author{%
  Jiahong Li\thanks{\raggedright Department of Applied Mathematics, The Hong Kong
  Polytechnic University, Hung Hom, Kowloon, Hong Kong.
  \mbox{E-mail: \href{mailto:jiahong24.li@polyu.edu.hk}{\nolinkurl{jiahong24.li@polyu.edu.hk}}}.}
  \and
  Defeng Sun\thanks{\raggedright Department of Applied Mathematics, The Hong Kong
  Polytechnic University, Hung Hom, Kowloon, Hong Kong.
  \mbox{E-mail: \href{mailto:defeng.sun@polyu.edu.hk}{\nolinkurl{defeng.sun@polyu.edu.hk}}}.}
  \and
  Hailing Wang\thanks{\raggedright Department of Applied Mathematics, The Hong Kong
  Polytechnic University, Hung Hom, Kowloon, Hong Kong.
  \mbox{E-mail: \href{mailto:hailing.wang@polyu.edu.hk}{\nolinkurl{hailing.wang@polyu.edu.hk}}}.}
  \and
  Liping Zhang\thanks{\raggedright Department of Mathematical Sciences, Tsinghua
  University, Beijing 100084, China. 
  \mbox{E-mail: \href{mailto:lipingzhang@tsinghua.edu.cn}{\nolinkurl{lipingzhang@tsinghua.edu.cn}}}.}
  \and
  Wei Zhao\thanks{\raggedright Department of Mathematical Sciences, Tsinghua University,
  Beijing 100084, China.
  \mbox{E-mail: \href{mailto:zhaow22@mails.tsinghua.edu.cn}{\nolinkurl{zhaow22@mails.tsinghua.edu.cn}}}.}%
}
\date{}

\newcommand{\keywords}[1]{\par\noindent\textbf{Key words.}\begingroup\def\and{\unskip, } #1\endgroup\par}
\newcommand{\subclass}[1]{\medskip\noindent\textbf{MSC codes.}\begingroup\def\and{\unskip, } #1\endgroup\par\medskip}

\begin{document}
\maketitle
\vspace{-1.5em}

\begin{abstract}
We introduce SHARP, a Structure-exploiting Halpern-Accelerated
Reflected splitting method with Preconditioning, for monotone
inclusions involving a monotone Lipschitz operator and a
maximal monotone operator. SHARP accommodates possibly degenerate
preconditioners and requires one forward evaluation and one
generalized backward step per iteration, without a separate
linear solve with the preconditioner. Under suitable step-size
and range conditions, SHARP achieves an $\mathcal{O}(1/k)$ bound on the last-iterate natural residual and converges in the induced seminorm.
Local continuity of the generalized backward map further ensures strong
convergence in the original Hilbert space. For affine forward operators in finite dimensions, a product-space
Halpern representation gives the same residual rate and strong
convergence under weaker step-size and range conditions. For structured problems, we
construct tailored preconditioners to simplify the generalized backward
step and improve convergence in practice. Numerical experiments on
large-scale structured problems demonstrate the computational efficiency
of the proposed method.

\end{abstract}

\keywords{Monotone inclusion \and Halpern acceleration \and
Degenerate preconditioning \and Operator splitting \and Structure-exploiting algorithms \and Large-scale structured programs}
\subclass{47H05 \and 49J40 \and 65K15 \and 90C25}

\section{Introduction}

In this paper, we consider the monotone inclusion problem in a real
Hilbert space $\mathcal H$:
\begin{equation}\label{eq:mi}
\text{find }x\in\mathcal H
\quad\text{such that}\quad
0\in F(x)+B(x),
\end{equation}
where $F:\mathcal H\to\mathcal H$ is monotone and Lipschitz continuous
with constant $L>0$, and
$B:\mathcal H\rightrightarrows\mathcal H$ is maximal monotone.
We assume throughout that the solution set
$\Scal:=(F+B)^{-1}(0)$ is nonempty.
Problem~\eqref{eq:mi} arises naturally from optimality conditions
for convex optimization and convex--concave saddle-point problems,
as well as from equilibrium conditions in monotone games~\cite{facchinei2003finite}.

For the special case $F=0$, preconditioned proximal-point
methods~\cite{rockafellar1976monotone,burke1999variable} use
problem-specific, possibly degenerate
preconditioners~\cite{bredies2022degenerate} to improve the scaling
of the iteration and simplify the backward step.
The resulting algorithms can be further accelerated by Halpern
iteration~\cite{halpern1967fixed,kim2021accelerated,sun2025accelerating},
which takes a convex combination of the initial point and the evaluation of the underlying nonexpansive operator.
With appropriate anchoring weights, Halpern iteration achieves both
explicit $\mathcal O(1/k)$ fixed-point residual
bounds~\cite{sabach2017first,lieder2021convergence} and strong
convergence~\cite{wittmann1992approximation}.
Combining Halpern acceleration with preconditioned proximal-point
methods has proven effective in practice, as demonstrated by HPR-LP~\cite{chen2026hprlp}, which has emerged as a leading GPU-based solver for large-scale Linear Programming (LP)~\cite{mittelmann2026lpfeas}.

Achieving accelerated convergence with structured preconditioning is more challenging when $F\neq0$. Applying an accelerated proximal-point method directly to $F+B$ generally requires evaluating the resolvent of their sum, which can be expensive even when the resolvent of $B$ is easy to compute. This difficulty motivates the use of splitting methods that treat $F$ explicitly while preserving a backward step involving $B$. However, retaining the computational advantages of preconditioning
in this setting depends on how the preconditioner enters
both the forward and backward steps.

For a bounded, self-adjoint, uniformly positive definite preconditioner $P$, the direct preconditioned formulation of the Modified Forward--Backward Splitting (MFBS) method of Tseng~\cite{tseng2000modified} takes the form
\[
y_k=(P+\gamma B)^{-1}\bigl(Px_k-\gamma F(x_k)\bigr),\qquad
x_{k+1}=y_k-\gamma P^{-1}\bigl(F(y_k)-F(x_k)\bigr).
\]
As in the case $F=0$, the preconditioner $P$ can be
chosen to simplify the generalized backward step.
However, the correction step requires a separate
application of $P^{-1}$, which may involve solving
a coupled linear system and offset the computational
benefits of preconditioning.
Moreover, this formulation does not directly
accommodate degenerate preconditioners.
These limitations motivate the use of the
Forward--Reflected--Backward (FRB) method of Malitsky and
Tam~\cite{malitsky2020forward}, whose direct preconditioned
iteration eliminates the separate correction step and takes the form
\[
x_{k+1}=(P+\gamma B)^{-1}\bigl(Px_k-2\gamma F(x_k)+\gamma F(x_{k-1})\bigr).
\]
Since $P$ appears only in $Px_k$ and $(P+\gamma B)^{-1}$,
the iteration remains well defined with a degenerate $P$
whenever the backward step has a unique solution.
A concrete example is the symmetric Gauss--Seidel (sGS) based alternating proximal augmented Lagrangian method of Sun, Wang, and Zhao~\cite{sun2026sgs}, which can be written as an FRB iteration with a degenerate preconditioner when the subproblems are solved exactly and the dual step-size is one.

While FRB accommodates flexible preconditioning, it does not achieve the accelerated $\mathcal O(1/k)$ last-iterate rate of Halpern-type methods\cite{golowich2020tight}. A representative splitting method that attains this rate is composite Fast Extra-Gradient (composite-FEG)~\cite{lee2021fast,cai2024accelerated}, which incorporates Halpern anchoring into an extra-gradient iteration~\cite{korpelevich1976extragradient}. For a uniformly positive definite preconditioner $P$, the direct
preconditioned iteration starts from an arbitrary $x_0$ with $d_0=0$
and computes, for $k\ge0$,
\[
y_k=\frac{1}{k+1}x_0+\frac{k}{k+1}x_k-\frac{\gamma k}{k+1}P^{-1}d_k,\qquad x_{k+1}=(P+\gamma B)^{-1}\!\left[P\left(\frac{1}{k+1}x_0+\frac{k}{k+1}x_k\right)-\gamma F(y_k)\right],
\]
where $d_{k+1}:=\frac{1}{\gamma}P\left(\frac{1}{k+1}x_0+\frac{k}{k+1}x_k-x_{k+1}\right)+F(x_{k+1})-F(y_k)\in(F+B)(x_{k+1})$.
Composite-FEG evaluates $F$ twice per iteration. Its prediction step
for $y_k$ requires a separate application of $P^{-1}$, limiting the
direct formulation to invertible preconditioners and potentially
offsetting the benefits of structured preconditioning.

In the monotone setting, the Accelerated Reflected Gradient (ARG) method of Cai and Zheng~\cite{cai2023accelerated} admits a direct preconditioned formulation that retains the same main update $x_{k+1}$ as composite-FEG, while updating the predictor via $y_0 = x_0$ and
\[
y_k=2x_k-x_{k-1}+\frac{1}{k+1}(x_0-x_k)-\frac{1}{k}(x_0-x_{k-1}),\qquad k\ge1.
\]
This affine prediction avoids a separate inverse-preconditioner operation and eliminates the need to evaluate $F(x_k)$ for the next prediction.
Each iteration therefore requires one forward evaluation and one generalized backward step, making ARG a promising candidate for accelerated splitting with low-cost, possibly degenerate preconditioning. However, for the identity preconditioner $P=\Id$, its original analysis establishes an $\mathcal O(1/k)$ last-iterate natural-residual bound only for $0<\gamma\le 1/(2\sqrt{6}L)$. It also does not establish convergence of the entire iterate sequence to a solution. A larger step-size range and convergence of the iterates are obtained by the Fast Reflected Forward--Backward (Fast RFB) method of Bo\c{t}, Nguyen, and Zong~\cite{bot2025fast}, which also admits a preconditioned formulation
with one forward evaluation and one generalized backward step per iteration. Specifically, given $z_0,u_0,w_0\in\mathcal H$, set
$z_1=(P+\gamma B)^{-1}\bigl(Pu_0-\gamma F(w_0)\bigr)$
and compute, for $k\ge1$,
\[
\begin{aligned}
u_k
&=z_k+\frac{k}{k+\alpha}(z_k-z_{k-1})
  +\left(1-\frac{c}{k+\alpha}\right)(u_{k-1}-z_k),\\
w_k
&=z_k+u_k-u_{k-1},\qquad
z_{k+1}=(P+\gamma B)^{-1}\bigl(Pu_k-\gamma F(w_k)\bigr).
\end{aligned}
\]
For $P=\Id$ and $0<\gamma<1/(2L)$, the original
analysis establishes weak convergence to a solution and asymptotic
$o(1/k)$ decay of the last-iterate natural residual under the
prescribed conditions on $\alpha$ and $c$. This asymptotic rate alone, however, does not yield an explicit bound on the
natural residual in terms of the iteration count and the step-size.

These developments motivate the search for a method that combines explicit last-iterate natural-residual bounds and strong convergence with a computational structure suitable for low-cost, possibly degenerate preconditioning.
To this end, we introduce \textbf{SHARP}, a \textbf{S}tructure-exploiting \textbf{H}alpern-\textbf{A}ccelerated \textbf{R}eflected splitting method with \textbf{P}reconditioning.
Each iteration evaluates $F$
once at a reflected predictor and performs one generalized backward
step with Halpern anchoring. The preconditioner $P$ enters only
through multiplication and the generalized backward map.

Our main contributions are as follows.
\begin{enumerate}
\renewcommand{\labelenumi}{(\roman{enumi})}

\item \textbf{An accelerated reflected splitting method with possibly
degenerate preconditioning.}
We introduce SHARP and, under mild conditions, establish an explicit
$\mathcal O(1/k)$ bound on the last-iterate natural residual and
convergence of the iterates in the seminorm induced by the fixed
preconditioner. Local continuity of the generalized backward map
further yields strong convergence in the original Hilbert space.
In a finite-dimensional space with $P=\Id$, SHARP can match the best-known
forward-oracle complexity bound for attaining a given tolerance
on the natural residual.

\item \textbf{A product-space Halpern representation for affine forward operators.} For affine forward operators in finite dimensions, we represent SHARP
as a Halpern iteration on a product space equipped with a seminorm.
This representation yields an explicit $\mathcal{O}(1/k)$ bound on the
last-iterate natural residual and strong convergence, including at the
step-size endpoint, under a relaxed range condition that allows $F$ to
take values outside $\operatorname{ran}P$.

\item \textbf{Structured preconditioners and numerical applications.} We exploit the block coupling in quadratic recourse games and the shared network structure in multicommodity-flow linear programs to construct preconditioners based on sGS decomposition and Schur-complement elimination, respectively. These designs simplify the generalized backward step through block-wise computation and reusable matrix factorizations. Numerical experiments on these problems and zero-sum network games demonstrate the efficiency of SHARP and the proposed preconditioners.

\end{enumerate}

The remainder of this paper is organized as follows.
Section~\ref{sec:arg} analyzes SHARP with the identity preconditioner. Section~\ref{sec:preconditioning} extends the analysis to possibly degenerate preconditioners through a range-space reduction. Section~\ref{sec:affine} develops the product-space analysis for affine forward operators. Sections~\ref{sec:preconditioner-design} and~\ref{sec:numerics} present structured preconditioners and numerical experiments, respectively. Section~\ref{sec:conclusion} concludes the paper.

\paragraph*{Notation.}
For a nonempty closed convex set $C$, $\Pi_C$ denotes the metric
projection onto $C$, and $\dist(x,C):=\norm{x-\Pi_C x}$ is the
distance from $x$ to $C$. We write $\mathcal N_C$ for the normal
cone operator of $C$. The $d\times d$ identity matrix is denoted
by $\mathrm{I}_d$. For $\gamma>0$,
$J_{\gamma B}:=(\Id+\gamma B)^{-1}$ denotes the resolvent of $\gamma B$.

\section{SHARP with the Identity Metric}\label{sec:arg}

We first introduce SHARP with the identity metric and establish an $\mathcal O(1/k)$ last-iterate natural-residual bound as well as strong convergence. These results provide the basis for studying the preconditioned method in Section~\ref{sec:preconditioning}. Throughout this section, we impose the following assumption. 

\begin{assumption}\label{ass:standing}
$F\colon\Hcal\to\Hcal$ is monotone and Lipschitz continuous with constant $L>0$, $B\colon\Hcal\rightrightarrows\Hcal$ is maximal monotone, and the solution set $\Scal$ of~\eqref{eq:mi} is nonempty.
\end{assumption}

SHARP evaluates $F$ at a reflected predictor and incorporates Halpern-type anchoring into the main update. The resulting iteration is presented in Algorithm~\ref{alg:arg}.

\begin{algorithm}[H]
\caption{SHARP with the Identity Metric for Solving~\eqref{eq:mi}}\label{alg:arg}
\begin{algorithmic}[1]
\STATE \textbf{Input:} Initial point $x_0\in\Hcal$ and step-size $\gamma>0$.
\STATE \textbf{Initialization:} Set $x_{-1}=x_0$ and $k=0$.
\STATE \textbf{Step 1:} Compute the reflected predictor:
\begin{equation}\label{eq:sharp-reflected-update}
 y_k=x_k+\frac{k-1}{k+1}(x_k-x_{k-1}).
\end{equation}
\STATE \textbf{Step 2:} Compute the main update:
\begin{equation}\label{eq:sharp-backward-update}
 x_{k+1}=J_{\gamma B}\left(\frac{1}{k+1}x_0+\frac{k}{k+1}x_k-\gamma F(y_k)\right).
\end{equation}
Set $k\leftarrow k+1$ and return to \textbf{Step 1}.
\STATE \textbf{Output:} The sequence $\{(x_k,y_k)\}_{k\ge0}$.
\end{algorithmic}
\end{algorithm}

Each iteration of Algorithm~\ref{alg:arg} requires one evaluation of $F$ and one application of $J_{\gamma B}$. Under Assumption~\ref{ass:standing}, the resolvent $J_{\gamma B}$ is everywhere defined, single-valued, and firmly nonexpansive for every $\gamma>0$~\cite{bauschke2017convex}, which ensures that the iteration is well defined.

To measure the last-iterate accuracy, we use the natural residual
$\rnat(x):=\norm{x-J_B(x-F(x))}$,
which satisfies $\rnat(x)=0$ if and only if $x\in \Scal$.
For convenience in the analysis, we introduce the following residual vector associated with the main update~\eqref{eq:sharp-backward-update}:
\begin{equation}\label{eq:sharp-dk}
d_{k+1}:=
\gamma^{-1}\left(
\frac{1}{k+1}x_0+\frac{k}{k+1}x_k-x_{k+1}
\right)
+F(x_{k+1})-F(y_k)
\in(F+B)(x_{k+1}),\qquad k\ge0.
\end{equation}
Since $J_B$ is nonexpansive and $x_k=J_B(x_k-F(x_k)+d_k)$, we have
$\rnat(x_k)\le\norm{d_k}$ for all $k\ge1$.
Our analysis relies on the Lyapunov inequality established in the following lemma.

\begin{lemma}\label{lem:sharp-energy}
Suppose that Assumption~\ref{ass:standing} holds.
For the sequences $\{(x_k,y_k)\}_{k\ge0}$ generated by Algorithm~\ref{alg:arg}
and $\{d_k\}_{k\ge1}$ defined as in~\eqref{eq:sharp-dk},
set $\mathcal L_0:=0$ and define
\begin{equation}\label{eq:sharp-energy}
 \mathcal L_k:=\frac{\gamma^2k^2}{2}
 \bigl(\norm{d_k}^2+\norm{F(x_k)-F(y_{k-1})}^2\bigr)
 +\gamma k\ip{d_k}{x_k-x_0},\qquad k\ge1.
\end{equation}
Then, for every $k\ge0$,
\begin{equation}\label{eq:sharp-descent}
 \mathcal L_k-\mathcal L_{k+1}
 \ge\left(\frac14-\gamma^2L^2\right)(k+1)^2\norm{x_{k+1}-y_k}^2.
\end{equation}
\end{lemma}

\begin{proof}
We first prove~\eqref{eq:sharp-descent} for $k=0$.
The initialization $x_{-1}=x_0$ and~\eqref{eq:sharp-reflected-update} give $y_0=x_0$, so definition~\eqref{eq:sharp-dk} yields
$\gamma d_1=x_0-x_1+\gamma(F(x_1)-F(x_0))$.
Substituting this identity into~\eqref{eq:sharp-energy} and expanding, we obtain
\[
\mathcal L_1=-\frac12\norm{x_1-x_0}^2
+\gamma^2\norm{F(x_1)-F(x_0)}^2.
\]
Hence, by the $L$-Lipschitz continuity of $F$ and $\mathcal L_0=0$, we have
\[
\mathcal L_0-\mathcal L_1
\ge\left(\frac12-\gamma^2L^2\right)\norm{x_1-y_0}^2
\ge\left(\frac14-\gamma^2L^2\right)\norm{x_1-y_0}^2,
\]
which establishes the initial case.

Now fix $k\ge1$. It follows from the identity~\eqref{eq:sharp-dk} that, for all $\ell\ge0$,
\begin{equation}\label{eq:sharp-proof-anchor}
\begin{aligned}
x_\ell-x_0
=-(\ell+1)(x_{\ell+1}-x_\ell)-\gamma(\ell+1)
\bigl(d_{\ell+1}-F(x_{\ell+1})+F(y_\ell)\bigr).
\end{aligned}
\end{equation}
Subtracting the identities for $\ell=k$ and $\ell=k-1$ yields
\begin{equation}\label{eq:sharp-proof-relation}
\begin{aligned}
&\gamma k\bigl(d_k-F(x_k)+F(y_{k-1})\bigr)
-\gamma(k+1)\bigl(d_{k+1}-F(x_{k+1})+F(y_k)\bigr)\\
&=(k+1)(x_{k+1}-x_k)-(k-1)(x_k-x_{k-1})=(k+1)(x_{k+1}-y_k),
\end{aligned}
\end{equation}
where the second equality follows
from~\eqref{eq:sharp-reflected-update}.

Expanding $\mathcal L_k-\mathcal L_{k+1}$
and using~\eqref{eq:sharp-proof-anchor}
to eliminate the inner-product terms involving $x_0$, we obtain
\begin{align*}
\mathcal L_k-\mathcal L_{k+1}
&=\gamma k(k+1)
\ip{d_{k+1}-d_k}{x_{k+1}-x_k}\\
&\quad+\frac{\gamma^2}{2}
\norm{
k d_k-(k+1)d_{k+1}
+(k+1)(F(x_{k+1})-F(y_k))
}^2\\
&\quad+\frac{\gamma^2k^2}{2}
\norm{F(x_k)-F(y_{k-1})}^2-\gamma^2(k+1)^2
\norm{F(x_{k+1})-F(y_k)}^2.
\end{align*}
By~\eqref{eq:sharp-proof-relation}, the second
term can be rewritten as $\frac12
\norm{
(k+1)(x_{k+1}-y_k)
+\gamma k(F(x_k)-F(y_{k-1}))
}^2$.
Completing the square then gives
\begin{align*}
\mathcal L_k-\mathcal L_{k+1}
&=\gamma k(k+1)
\ip{d_{k+1}-d_k}{x_{k+1}-x_k}\\
&\quad+\frac14
\norm{
(k+1)(x_{k+1}-y_k)
+2\gamma k(F(x_k)-F(y_{k-1}))
}^2\\
&\quad+\frac{(k+1)^2}{4}
\norm{x_{k+1}-y_k}^2-\gamma^2(k+1)^2
\norm{F(x_{k+1})-F(y_k)}^2.
\end{align*}
By~\eqref{eq:sharp-dk} and the monotonicity of $F+B$,
the first term is nonnegative. Dropping the first two terms
and using the $L$-Lipschitz continuity of $F$, we conclude that
\[
\mathcal L_k-\mathcal L_{k+1}
\ge
\left(\frac14-\gamma^2L^2\right)
(k+1)^2\norm{x_{k+1}-y_k}^2,
\]
which proves~\eqref{eq:sharp-descent}.
\end{proof}

We next establish an $\mathcal{O}(1/k)$ last-iterate natural-residual bound and strong convergence of the iterates. Inspired by the merging path property of Yoon and Ryu~\cite{yoon2025flock}, we compare the SHARP iterates with a reference sequence converging to $\Pi_{\Scal}x_0$ and show that the distance between the two sequences tends to zero.

\begin{theorem}\label{thm:sharp-convergence}
Suppose that Assumption~\ref{ass:standing} holds. Then the sequences $\{(x_k,y_k)\}_{k\ge0}$ generated by Algorithm~\ref{alg:arg}
and $\{d_k\}_{k\ge1}$ defined as in~\eqref{eq:sharp-dk} satisfy
the following properties.
\begin{enumerate}
  \renewcommand{\labelenumi}{(\roman{enumi})}
  \item If $0<\gamma\le1/(2L)$, it holds that
  \begin{equation}\label{eq:sharp-rate}
   \rnat(x_k)\le\norm{d_k}\le\frac{2\dist(x_0,\Scal)}{\gamma k},\qquad k\ge1.
  \end{equation}
  \item If $0<\gamma<1/(2L)$, then
  $x_k,y_k\to\bar x:=\Pi_{\Scal}x_0$ strongly in $\Hcal$.
\end{enumerate}
\end{theorem}

\begin{proof}
\noindent(i) Suppose that $0<\gamma\le1/(2L)$. By
Lemma~\ref{lem:sharp-energy}, $\mathcal L_k\le\mathcal L_0=0$
for every $k\ge1$. For any $s\in\Scal$, since $d_k\in(F+B)(x_k)$ and
$0\in(F+B)(s)$, the monotonicity of $F+B$
gives $\ip{d_k}{x_k-s}\ge0$. Combining these inequalities with
\eqref{eq:sharp-energy} yields
\begin{equation}\label{eq:sharp-residual-estimate}
\frac{\gamma k}{2}\norm{d_k}^2
\le \ip{d_k}{x_0-x_k}
\le \ip{d_k}{x_0-s},
\qquad s\in\Scal,\quad k\ge1.
\end{equation}
Under Assumption~\ref{ass:standing}, $F+B$ is maximal monotone
and its zero set $\Scal$ is closed and convex~\cite{bauschke2017convex}.
Set $\bar x:=\Pi_{\Scal}x_0$. Taking $s=\bar x$ in
\eqref{eq:sharp-residual-estimate} and applying Cauchy--Schwarz gives
$\gamma k\norm{d_k}\le2\norm{x_0-\bar x}=2\dist(x_0,\Scal)$.
Together with $\rnat(x_k)\le\norm{d_k}$, this proves
\eqref{eq:sharp-rate}.

\noindent(ii) Suppose now that $0<\gamma<1/(2L)$. We first establish a uniform lower bound for $\mathcal L_k$. Using the equality~\eqref{eq:sharp-energy} and the monotonicity inequality
$\ip{d_k}{x_k-\bar x}\ge0$, we obtain
\[
 \mathcal L_k
 \ge
 \frac{\gamma^2k^2}{2}\norm{d_k}^2
 +\gamma k\ip{d_k}{\bar x-x_0}
 \ge-\frac{\norm{x_0-\bar x}^2}{2}.
\]
Summing~\eqref{eq:sharp-descent} from $k=0$ to $k=N-1$, we obtain
\[
\left(\frac14-\gamma^2L^2\right)
\sum_{k=0}^{N-1}(k+1)^2\norm{x_{k+1}-y_k}^2
\le\mathcal L_0-\mathcal L_N
\le\frac12\norm{x_0-\bar x}^2.
\]
Since $\gamma<1/(2L)$, letting $N\to\infty$ yields $\sum_{k=0}^\infty(k+1)^2\norm{x_{k+1}-y_k}^2<\infty$. 
Consequently, we have $(k+1)\norm{x_{k+1}-y_k}\to0$.

We define the reference sequence $\{\tilde x_k\}_{k\ge0}$ by
\begin{equation}\label{eq: reference-iterate-update}
\tilde x_0=x_0,\qquad
\tilde x_{k+1}=J_{\gamma(F+B)}\left(
\frac{1}{k+1}x_0+\frac{k}{k+1}\tilde x_k\right),\qquad k\ge0.
\end{equation}
Set $a_k:=(x_0+k\tilde x_k)/(k+1)$.
Then $a_0=x_0$, $\tilde x_{k+1}=J_{\gamma(F+B)}(a_k)$, and $a_{k+1}=\frac{1}{k+2}a_0
+\frac{k+1}{k+2}J_{\gamma(F+B)}(a_k)$, for all $k\ge0$.
Thus, $\{a_k\}_{k \ge 0}$ is a Halpern iteration for the nonexpansive map $J_{\gamma(F+B)}$,
whose fixed-point set is $\Scal$. By Wittmann's theorem~\cite[Theorem~2]{wittmann1992approximation}, this sequence converges strongly to $\bar x$. Combining this with the nonexpansiveness of $J_{\gamma(F+B)}$, we then obtain
$\tilde x_k\to\bar x$ strongly.

It remains to show that $\norm{x_k-\tilde x_k}\to0$. For all $k\ge0$, the backward update~\eqref{eq:sharp-backward-update} can be rewritten as
\[
x_{k+1}=J_{\gamma(F+B)}\left(
\frac{1}{k+1}x_0+\frac{k}{k+1}x_k
+\gamma\bigl(F(x_{k+1})-F(y_k)\bigr)\right).
\]
Comparing the above identity with the reference iteration~\eqref{eq: reference-iterate-update} and using the nonexpansiveness of $J_{\gamma(F+B)}$ and the $L$-Lipschitz continuity of $F$, we obtain
\[
(k+1)\norm{x_{k+1}-\tilde x_{k+1}}
\le k\norm{x_k-\tilde x_k}
+\gamma L(k+1)\norm{x_{k+1}-y_k}.
\]
Summing these inequalities from $k=0$ to $k=N-1$ and dividing
by $N$, we have
\[
\norm{x_N-\tilde x_N}
\le\frac{\gamma L}{N}
\sum_{k=0}^{N-1}(k+1)\norm{x_{k+1}-y_k},
\qquad N\ge1.
\]
Since $(k+1)\norm{x_{k+1}-y_k}\to0$, the right-hand side tends
to zero, which implies $\norm{x_N-\tilde x_N}\to0$.
Together with $\tilde x_k\to\bar x$, we conclude that $x_k$ converges strongly to $\bar x$ in $\Hcal$. The reflected update~\eqref{eq:sharp-reflected-update} then gives
$y_k\to\bar x$ strongly in $\Hcal$.
\end{proof}

\begin{remark}[Analysis of forward-oracle complexity]
\label{rem:forward-oracle-complexity}
In a finite-dimensional space with $P=\Id$, let $\{u_k\}_{k \ge 0}$ denote
the composite-FEG iterates initialized at $u_0=x_0$.
The analysis in~\cite[Section~4.1]{cai2024accelerated} yields
the last-iterate natural residual bound
\[
\left\|u_k-J_B\bigl(u_k-F(u_k)\bigr)\right\|
\le \frac{2\,\dist(x_0,\Scal)}{\gamma k},
\qquad k\ge1,\quad 0<\gamma\le\frac{1}{L}.
\]
With $\gamma=1/L$ and two forward evaluations per iteration,
composite-FEG reaches a prescribed tolerance $\varepsilon>0$
on the natural residual within
$2\bigl\lceil 2L\,\dist(x_0,\Scal)/\varepsilon\bigr\rceil$
forward evaluations. To the best of our knowledge, composite-FEG
gives the lowest bound among prior first-order splitting methods
for this class of monotone inclusions. For SHARP, applying~\eqref{eq:sharp-rate} with $\gamma=1/(2L)$
and counting one forward evaluation per iteration gives the bound
$\bigl\lceil 4L\,\dist(x_0,\Scal)/\varepsilon\bigr\rceil$
for reaching the same tolerance, which is no larger than the
corresponding bound for composite-FEG.
\end{remark}

\begin{remark}[Relation to anchored FRBS]\label{rem:anchored-frbs}
For $k\ge1$, the updates in Algorithm~\ref{alg:arg} fit the anchored FRBS framework of Tran-Dinh and Nguyen-Trung~\cite[Equation~(22)]{trandinh2025accelerated}, with $\tau_k=1/(k+1)$, $\widehat\eta_k=\gamma(1-\tau_k)$, and $\beta_k=0$. However, the convergence result in \cite[Corollary~4.2]{trandinh2025accelerated} assumes $\tau_k=1/(k+\nu)$ with $\nu>2$ and a step-size of at most $1/(2\sqrt{2}L)$, so it does not directly apply to Algorithm~\ref{alg:arg} with $\gamma\in(0,1/(2L)]$.
\end{remark}

\section{SHARP with Possibly Degenerate Preconditioners}\label{sec:preconditioning}

We now study SHARP with a fixed, possibly degenerate
preconditioner. In addition to Assumption~\ref{ass:standing}, we impose the following conditions
to ensure that the generalized backward step is well defined.

\begin{assumption}\label{ass:degenerate}
Let $P\colon\Hcal\to\Hcal$ be a nonzero, bounded, self-adjoint, positive semidefinite linear operator with closed range. For the chosen
step-size $\gamma>0$, we assume that
\begin{enumerate}
\renewcommand{\labelenumi}{(\roman{enumi})}
\item $\ran F\subseteq\ran P$.
\item $(P+\gamma B)^{-1}$ is everywhere defined and single-valued on $\ran P$.
\end{enumerate}
\end{assumption}

Assumption~\ref{ass:degenerate}(ii) is equivalent to requiring that
$P$ be admissible for $\gamma B$ in the sense of
\cite[Definition~2.1]{bredies2022degenerate}. The iteration of SHARP is presented in Algorithm~\ref{alg:darg}. By Assumption~\ref{ass:degenerate}(i), the input to the backward
step~\eqref{eq:darg-backward-update} belongs to $\ran P$.
Assumption~\ref{ass:degenerate}(ii) then ensures that $x_{k+1}$ is uniquely determined.

\begin{algorithm}[H]
\caption{SHARP with Degenerate Preconditioning for Solving~\eqref{eq:mi}}\label{alg:darg}
\begin{algorithmic}[1]
\STATE \textbf{Input:} Initial point $x_0\in\Hcal$, step-size $\gamma>0$,
and preconditioner $P$.
\STATE \textbf{Initialization:} Set $x_{-1}=x_0$ and $k=0$. 
\STATE \textbf{Step 1:} Compute the reflected predictor:
\begin{equation}\label{eq:darg-reflected-update}
 y_k=x_k+\frac{k-1}{k+1}(x_k-x_{k-1}).
\end{equation}
\STATE \textbf{Step 2:} Compute the preconditioned main update:
\begin{equation}\label{eq:darg-backward-update}
 x_{k+1}=(P+\gamma B)^{-1}\left(
 P\left(\frac{1}{k+1}x_0+\frac{k}{k+1}x_k\right)
 -\gamma F(y_k)\right).
\end{equation}
Set $k\leftarrow k+1$ and return to \textbf{Step 1}.
\STATE \textbf{Output:} The sequence $\{(x_k,y_k)\}_{k\ge0}$.
\end{algorithmic}
\end{algorithm}

As in the identity-metric setting, rearranging the update~\eqref{eq:darg-backward-update} yields the residual vector
\begin{equation}\label{eq:darg-dk}
d_{k+1}:=
\gamma^{-1}P\left(
\frac{1}{k+1}x_0+\frac{k}{k+1}x_k-x_{k+1}
\right)
+F(x_{k+1})-F(y_k)
\in(F+B)(x_{k+1}),\qquad k\ge0.
\end{equation}

Following the reduced-space approach of~\cite{bredies2022degenerate},
we work on the closed subspace $\Vcal:=\ran P=(\ker P)^\perp$.
Since $P$ has closed range by Assumption~\ref{ass:degenerate}, the restriction $P_{\Vcal}:=P|_{\Vcal}$ is self-adjoint, uniformly
positive, and boundedly invertible on $\Vcal$. We equip $\Vcal$ with
the inner product $\ip{z}{w}_{\Vcal}:=\ip{P_{\Vcal}^{-1}z}{w}$.
The induced norm $\norm{z}_{\Vcal}=\norm{P_{\Vcal}^{-1/2}z}$ satisfies
\begin{equation}\label{eq:range-norm-equivalence}
\frac{\norm{z}}{\sqrt{\norm{P}}}
\le \norm{z}_{\Vcal}
\le \sqrt{\norm{P_{\Vcal}^{-1}}}\,\norm{z},
\qquad z\in\Vcal.
\end{equation}
These bounds and the closedness of $\Vcal$ imply that
$(\Vcal,\ip{\cdot}{\cdot}_{\Vcal})$ is a Hilbert space.
Write $\norm{h}_P:=\sqrt{\ip{Ph}{h}}$ for the seminorm induced by $P$.
By self-adjointness of $P$, for $b\in\Vcal$ and $h\in\Hcal$,
\begin{equation}\label{eq:range-metric-identities}
\ip{b}{Ph}_{\Vcal}
=\ip{PP_{\Vcal}^{-1}b}{h}
=\ip{b}{h},
\qquad
\norm{Ph}_{\Vcal}=\norm{h}_P.
\end{equation}

The following lemma relates the original inclusion to its reduction
on $(\Vcal,\ip{\cdot}{\cdot}_{\Vcal})$.

\begin{lemma}
\label{lem:exact-reduction}
Suppose that Assumptions~\ref{ass:standing}
and~\ref{ass:degenerate} hold. Then the following statements hold:
\begin{enumerate}
\renewcommand{\labelenumi}{(\roman{enumi})}

\item Let $F_{\Vcal}\colon\Vcal\to\Vcal$ be defined by
$F_{\Vcal}:=F\circ P_{\Vcal}^{-1}$.
Then $F_{\Vcal}(Px)=F(x)$ for every $x\in\Hcal$.
Moreover, $F_{\Vcal}$ is monotone and admits a Lipschitz constant
$L_P$ with $0<L_P\le L\norm{P_{\Vcal}^{-1}}$.

\item Define $B_{\Vcal}\colon\Vcal\rightrightarrows\Vcal$ by
$B_{\Vcal}(z):=(B\circ P^{-1})(z)\cap\Vcal$, where $P^{-1}$
denotes the set-valued inverse of $P$.
Set $\mathcal R_\gamma:=(P+\gamma B)^{-1}|_{\Vcal}$.
Then $B_{\Vcal}$ is maximal monotone and
$J_{\gamma B_{\Vcal}}=P\mathcal R_\gamma$.

\item The reduced solution set
$\Scal_{\Vcal}:=(F_{\Vcal}+B_{\Vcal})^{-1}(0)=P\Scal$
is nonempty, closed, and convex. The map $s\mapsto Ps$ is a
bijection from $\Scal$ onto $\Scal_{\Vcal}$ with inverse
$z\mapsto\mathcal R_\gamma(z-\gamma F_{\Vcal}(z))$.
\end{enumerate}
\end{lemma}

\begin{proof}
\noindent(i)
We first show that $F(x)=F(y)$ whenever $Px=Py$.
For any $h\in\Hcal$ and $t>0$, the monotonicity of $F$ gives
\[
0\le\ip{F(y+th)-F(x)}{y+th-x}
=t\ip{F(y+th)-F(x)}{h},
\]
where the equality uses $y-x\in\ker P=\Vcal^\perp$ and $\ran F\subseteq \Vcal$.
Dividing by $t$ and letting $t\downarrow0$, we obtain
$\ip{F(y)-F(x)}{h}\ge0$ for all $h\in\Hcal$, which implies $F(y)=F(x)$.
Taking $y=P_{\Vcal}^{-1}(Px)$ gives $Py=Px$ and hence
$F_{\Vcal}(Px)=F(y)=F(x)$.

For any $z=Px$ and $w=Py$ in $\Vcal$,~\eqref{eq:range-metric-identities}
and the monotonicity of $F$ give
\[
\ip{F_{\Vcal}(z)-F_{\Vcal}(w)}{z-w}_{\Vcal}
=\ip{F(x)-F(y)}{x-y}\ge0.
\]
To derive the Lipschitz bound, choose the
representatives $x=P_{\Vcal}^{-1}z$ and $y=P_{\Vcal}^{-1}w$ in $\Vcal$. Using
$\norm{P_{\Vcal}^{-1}b}=\norm{P_{\Vcal}^{-1/2}P_{\Vcal}^{-1/2}b}
\le\sqrt{\norm{P_{\Vcal}^{-1}}}\,\norm{b}_{\Vcal}$ for $b=z-w\in \Vcal$ together
with~\eqref{eq:range-norm-equivalence}, we obtain
\[
\norm{F_{\Vcal}(z)-F_{\Vcal}(w)}_{\Vcal}
\le\sqrt{\norm{P_{\Vcal}^{-1}}}\,L\norm{x-y}
\le L\norm{P_{\Vcal}^{-1}}\norm{z-w}_{\Vcal}.
\]

\noindent(ii)
Take $b\in B_{\Vcal}(z)$ and $c\in B_{\Vcal}(w)$, with representatives
$x,y\in\Hcal$ satisfying $Px=z$, $Py=w$, $b\in B(x)$, and $c\in B(y)$.
By~\eqref{eq:range-metric-identities} and the monotonicity of $B$, we have $\ip{b-c}{z-w}_{\Vcal}=\ip{b-c}{x-y}\ge0$. It follows that $B_{\Vcal}$ is monotone in $(\Vcal,\ip{\cdot}{\cdot}_{\Vcal})$.
For any $q\in \Vcal$, set $p=\mathcal R_\gamma(q)$.
Since $q,Pp\in \Vcal$, the inclusion $q\in(P+\gamma B)(p)$ implies
$(q-Pp)/\gamma\in B(p)\cap \Vcal$, and hence
\[
q\in(\Id_{\Vcal}+\gamma B_{\Vcal})(Pp).
\]
Thus $\ran(\Id_{\Vcal}+\gamma B_{\Vcal})=\Vcal$. Minty's theorem~\cite{bauschke2017convex}
therefore yields the maximal monotonicity of $B_{\Vcal}$ and the resolvent
identity $J_{\gamma B_{\Vcal}}=P\mathcal R_\gamma$.

\noindent(iii)
By part~(i) and (ii), $z\in \Scal_{\Vcal}$ if and only if there
exists $s\in\Hcal$ such that $Ps=z$ and $-F_{\Vcal}(z)=-F(s)\in B(s)$,
hence $\Scal_{\Vcal}=P\Scal$. Every such $s\in \Scal$ satisfies
\[
z-\gamma F_{\Vcal}(z)=Ps-\gamma F(s)\in(P+\gamma B)(s),
\]
and Assumption~\ref{ass:degenerate}(ii) uniquely determines
$s=\mathcal R_\gamma(z-\gamma F_{\Vcal}(z))$. Finally, since $F_{\Vcal}$ is
continuous, monotone, and everywhere defined on $\Vcal$ and $B_{\Vcal}$ is
maximal monotone, $F_{\Vcal}+B_{\Vcal}$ is maximal monotone~\cite{bauschke2017convex}.
Its zero set $\Scal_{\Vcal}$ is therefore closed and convex, and it is nonempty
because $\Scal\ne\varnothing$.
\end{proof}

We now write the iteration in $\Vcal$. For the iterates of Algorithm~\ref{alg:darg}, set $z_k:=Px_k$
and $v_k:=Py_k$, with $z_{-1}=z_0$. By
Lemma~\ref{lem:exact-reduction}, the reduced updates for $k\ge0$ are
\begin{equation}\label{eq:darg-reduced-iteration}
v_k = z_k+\frac{k-1}{k+1}(z_k-z_{k-1}),\qquad q_k = \frac{z_0+kz_k}{k+1}-\gamma F_{\Vcal}(v_k),
\qquad z_{k+1}=J_{\gamma B_{\Vcal}}(q_k).
\end{equation}
The original iterate is recovered as $x_{k+1}=\mathcal R_\gamma(q_k)$.
The pairs $(z_k,v_k)$ follow the identity-metric SHARP iteration for
$0\in(F_{\Vcal}+B_{\Vcal})(z)$ in $(\Vcal,\ip{\cdot}{\cdot}_{\Vcal})$.
Moreover, its residual vector coincides with the vector $d_{k+1}$
defined in~\eqref{eq:darg-dk}:
\begin{equation*}\label{eq:range-residual}
d_{k+1}
=\gamma^{-1}(q_k-z_{k+1})+F_{\Vcal}(z_{k+1})
\in(F_{\Vcal}+B_{\Vcal})(z_{k+1}),\qquad k\ge0.
\end{equation*}
In particular, $\norm{d_{k+1}}_{\Vcal}$ is well defined.

We now establish the main result of this section by applying the
identity-metric analysis to the reduced iteration on $\Vcal$.

\begin{theorem}\label{thm:darg}
Suppose that Assumptions~\ref{ass:standing}
and~\ref{ass:degenerate} hold. Let $L_P>0$ be a Lipschitz constant
of $F_{\Vcal}$ in $\norm{\cdot}_{\Vcal}$, and let $\{(x_k,y_k)\}_{k\ge0}$ be
generated by Algorithm~\ref{alg:darg} from $x_0\in\Hcal$. Set
$\rho_{P,0}:=\inf_{s\in \Scal}\norm{x_0-s}_P$. Then the following
statements hold:
\begin{enumerate}
\renewcommand{\labelenumi}{(\roman{enumi})}
\setlength{\itemsep}{0pt}
\setlength{\parsep}{0pt}

\item
If $0<\gamma\le1/(2L_P)$, then, for every $k\ge1$,
$\norm{d_k}_{\Vcal}\le 2\rho_{P,0}/(\gamma k)$ and
\begin{equation}\label{eq:darg-rate}
\rnat(x_k)
\le\frac{1}{\gamma k}\inf_{s\in \Scal}
\left\{
\norm{P(x_0-s)}
+\sqrt{\norm{P}}\,\norm{x_0-s}_P
\right\}
\le\frac{2\sqrt{\norm{P}}\,\rho_{P,0}}{\gamma k}.
\end{equation}

\item If $0<\gamma<1/(2L_P)$, then $x_k$ and $y_k$ converge in
$\norm{\cdot}_P$ to the unique minimizer $\bar x$ of
$\norm{x-x_0}_P^2$ over $x\in\Scal$.

\item Under the conditions of part~(ii), if $\mathcal R_\gamma$ is continuous at
$P\bar x-\gamma F(\bar x)$ in the ambient norm, then $x_k$ and $y_k$
converge strongly to $\bar x$ in $\Hcal$.

\end{enumerate}
\end{theorem}

\begin{proof}

By Lemma~\ref{lem:exact-reduction}
and~\eqref{eq:darg-reduced-iteration}, Theorem~\ref{thm:sharp-convergence} applies to $\{(z_k,v_k)\}_{k\ge0}$ in $(\Vcal,\ip{\cdot}{\cdot}_{\Vcal})$.
Moreover,~\eqref{eq:range-metric-identities} and $\Scal_{\Vcal}=P\Scal$ give
\[
\inf_{z\in \Scal_{\Vcal}}\norm{Px_0-z}_{\Vcal}
=\inf_{s\in \Scal}\norm{x_0-s}_P=\rho_{P,0}.
\]

\noindent(i) 
Theorem~\ref{thm:sharp-convergence}(i) gives
$\norm{d_k}_{\Vcal}\le 2\rho_{P,0}/(\gamma k)$.
For any $s\in\Scal$, applying~\eqref{eq:sharp-residual-estimate}
to the reduced iteration with $Ps\in\Scal_{\Vcal}$ and multiplying
by $\gamma k$ gives
\[
\frac{\gamma^2k^2}{2}\norm{d_k}_{\Vcal}^2
-\gamma k\ip{d_k}{P(x_0-s)}_{\Vcal}\le0,
\qquad s\in \Scal.
\]
Completing the square and taking square roots, we arrive at
\[
\norm{\gamma k d_k-P(x_0-s)}_{\Vcal}
\le\norm{x_0-s}_P,
\qquad s\in \Scal.
\]
Combining this with~\eqref{eq:range-norm-equivalence} and the triangle inequality, we have
\[
\gamma k\norm{d_k}
\le\norm{P(x_0-s)}
+\sqrt{\norm{P}}\,\norm{x_0-s}_P
\le 2\sqrt{\norm{P}}\,\norm{x_0-s}_P.
\]
Finally, taking the infimum over $s\in \Scal$ and utilizing $\rnat(x_k)\le\norm{d_k}$ gives the bounds~\eqref{eq:darg-rate}.

\noindent(ii)
Let $\bar z$ be the projection of $z_0$ onto $\Scal_{\Vcal}$ with
respect to $\norm{\cdot}_{\Vcal}$. By
Theorem~\ref{thm:sharp-convergence}(ii), $z_k$ and $v_k$ converge to
$\bar z$ in $\norm{\cdot}_{\Vcal}$. Lemma~\ref{lem:exact-reduction}(iii) gives a
unique $\bar x\in\Scal$ such that $P\bar x=\bar z$. Since
$\norm{x-x_0}_P^2=\norm{Px-z_0}_{\Vcal}^2$
by~\eqref{eq:range-metric-identities}, $\bar x$ is the unique minimizer
of $\norm{x-x_0}_P^2$ over $x\in\Scal$. The same norm identity gives
\[
\norm{x_k-\bar x}_P=\norm{z_k-\bar z}_{\Vcal}\longrightarrow0,
\qquad
\norm{y_k-\bar x}_P=\norm{v_k-\bar z}_{\Vcal}\longrightarrow0.
\]

\noindent(iii)
Part~(ii) gives $z_k,v_k\to\bar z=P\bar x$ in
$\norm{\cdot}_{\Vcal}$. The definition of $q_k$
in~\eqref{eq:darg-reduced-iteration} and the continuity of
$F_{\Vcal}$ yield $q_k\to\bar q$ in the same norm, where $
\bar q:=\bar z-\gamma F_{\Vcal}(\bar z)
       =P\bar x-\gamma F(\bar x).
$
By~\eqref{eq:range-norm-equivalence}, $q_k\to\bar q$ also holds
in the ambient norm. Continuity of $\mathcal R_\gamma$ at $\bar q$
and Lemma~\ref{lem:exact-reduction}(iii) give
\[
x_{k+1}=\mathcal R_\gamma(q_k)\longrightarrow
\mathcal R_\gamma(\bar q)=\bar x
\qquad\text{in }\Hcal.
\]
The reflected update~\eqref{eq:darg-reflected-update}
then yields $y_k\to\bar x$ in $\Hcal$.
\end{proof}

\section{SHARP with Affine Forward Operators}\label{sec:affine}

In this section, we consider the structured setting where $\Hcal=\R^n$ and $F$ is affine.
In this special case, we relax
the range requirement $\ran F\subseteq\ran P$ and establish an explicit $\mathcal O(1/k)$ last-iterate natural-residual bound and strong convergence, with equality allowed in the step-size condition. Throughout this section, we impose the following assumption.

\begin{assumption}\label{ass:affine}
Suppose that Assumption~\ref{ass:standing} holds in $\Hcal=\R^n$ with an affine operator $F(x)=Ax+c$, where $A\in\R^{n\times n}$ and $c\in\R^n$. Let $P\in\R^{n\times n}$ be a symmetric positive semidefinite matrix. For the chosen step-size $\gamma>0$, we assume that there exists a symmetric positive semidefinite matrix $\Sigma\in\R^{n\times n}$ such that
$ P_{\Sigma}:=P+\frac{\gamma}{2}\Sigma $ is nonzero and the following conditions hold:
\begin{enumerate}
\renewcommand{\labelenumi}{(\roman{enumi})}
\setlength{\itemsep}{0pt}
\setlength{\parsep}{0pt}
\item $F+B-\Sigma$ is monotone.
\item $(P+\gamma B)^{-1}$ is everywhere defined and single-valued on $\ran P_{\Sigma}-\{\gamma c\}$.
\end{enumerate}
\end{assumption}

From Assumption~\ref{ass:affine}(i), we notice that $\Sigma$ provides a monotonicity lower bound for $F+B$, while $P_{\Sigma}$ combines this bound with the preconditioner $P$. The backward input in Algorithm~\ref{alg:darg} has the form $Pu-\gamma Aw-\gamma c$.
By Assumption~\ref{ass:affine}(ii), the backward step is uniquely
solvable whenever $Pu-\gamma Aw\in\ran P_{\Sigma}$. To guarantee this range inclusion and establish nonexpansiveness of the product-space map, we impose the step-size condition
\begin{equation}\label{eq:affine-bilinear}
2\gamma\lvert\ip{Au}{v}\rvert
\le\norm{u}_{P_{\Sigma}}\norm{v}_{P_{\Sigma}},
\qquad u,v\in\R^n, 
\end{equation}
where $\norm{x}_{P_{\Sigma}}:=\sqrt{\ip{x}{P_{\Sigma}x}}$ denotes the seminorm induced by $P_{\Sigma}$. 

\subsection{Halpern Representation and Nonexpansiveness}

To reveal the hidden Halpern structure of Algorithm~\ref{alg:darg}
when $F$ is affine, we introduce the product-space mapping
$\mathcal T_\gamma\colon\R^n\times\R^n\to\R^n\times\R^n$ by
\begin{equation}\label{eq:affine-product-map}
\mathcal T_\gamma(u,v)
:=
\left(
(P+\gamma B)^{-1}\bigl(Pu-\gamma F(2u-v)\bigr),
\,u
\right).
\end{equation}
Then the iteration of the preconditioned reflected forward--backward
method~\cite{cevher2021reflected} satisfies
$(w_{k+1},w_k)=\mathcal T_\gamma(w_k,w_{k-1})$ with $w_{-1}=w_0$.
The identity $F(2u-v)=2F(u)-F(v)$ shows that $\mathcal T_\gamma$
also gives the preconditioned forward--reflected--backward
update~\cite{malitsky2020forward}. The following lemma establishes well-posedness of
Algorithm~\ref{alg:darg} and a Halpern representation
of an auxiliary sequence generated by SHARP. 

\begin{lemma}\label{lem:affine-well-posedness}
Suppose that Assumption~\ref{ass:affine} and~\eqref{eq:affine-bilinear} hold. Then the following statements hold:
\begin{enumerate}
\renewcommand{\labelenumi}{(\roman{enumi})}
\setlength{\itemsep}{0pt}
\setlength{\parsep}{0pt}
\item
The map $\mathcal T_\gamma$ defined
in~\eqref{eq:affine-product-map} is everywhere defined and single-valued on $\R^n\times\R^n$, and the update formulas of Algorithm~\ref{alg:darg} generate a unique sequence
$\{(x_k,y_k)\}_{k\ge0}$ for every initial point $x_0\in\R^n$.
\item
Let $\{(x_k,y_k)\}_{k\ge0}$ be generated by Algorithm~\ref{alg:darg}, and define auxiliary sequences
\begin{equation*}\label{eq:affine-auxiliary-points}
u_k:=\frac{1}{k+1}x_0 + \frac{k}{k+1}x_k,
\qquad
v_k:=2u_k-y_k,
\qquad k\ge0.
\end{equation*}
Then $(u_0,v_0)=(x_0,x_0)$, and the sequence
$\{(u_k,v_k)\}_{k\ge0}$ satisfies the following Halpern recursion:
\begin{equation}\label{eq:affine-halpern}
(u_{k+1},v_{k+1})
=\frac{1}{k+2}(x_0,x_0)
+\frac{k+1}{k+2}\mathcal T_\gamma(u_k,v_k),\qquad k\geq 0.
\end{equation}
\end{enumerate}
\end{lemma}

\begin{proof}
\noindent(i)
Noting that $P,\Sigma\succeq0$, we have
\begin{equation}\label{inclu: KP}
\ker P_{\Sigma}=\ker P\cap\ker\Sigma,
\qquad
\ran P\subseteq\ran P_{\Sigma}.
\end{equation}
For $z\in\ker P_{\Sigma}$, condition~\eqref{eq:affine-bilinear}
gives $\ip{Az}{u}=\ip{Au}{z}=0$ for every $u\in\R^n$.
Hence,
\begin{equation}\label{eq: AZU=0}
Az=A^\top z=0,
\qquad z\in\ker P_{\Sigma},
\end{equation}
and therefore $\ran A\subseteq(\ker P_{\Sigma})^\perp=\ran P_{\Sigma}$.
It follows that
\[
Pu-\gamma F(w)\in\ran P_{\Sigma}-\{\gamma c\},
\qquad u,w\in\R^n.
\]
Assumption~\ref{ass:affine}(ii) now establishes
the claimed well-posedness.

\noindent(ii)
The initialization gives $u_0=v_0=x_0$.
Since $y_k=2u_k-v_k$, the backward update yields
\begin{equation}\label{eq:affine-output}
\mathcal T_\gamma(u_k,v_k)=(x_{k+1},u_k),
\qquad k\ge0.
\end{equation}
The reflected update also gives
\[
v_{k+1} = \frac{2}{k+2}x_0 + \frac{k}{k+2} x_k = \frac{1}{k+2}x_0 + \frac{k+1}{k+2}u_k
\]
Combining this identity with the definition of
$u_{k+1}$ proves~\eqref{eq:affine-halpern}.
\end{proof}

The Halpern representation~\eqref{eq:affine-halpern} remains valid for
general $F$ whenever the backward steps are well defined.
For affine $F(x)=Ax+c$, define the symmetric matrix $M$ by
\begin{equation}\label{eq:affine-product-metric}
\ip{(a,b)}{M(a,b)}
=\norm{a}_{P+\gamma\Sigma}^2
-2\gamma\ip{a}{A(a-b)}
+\frac12\norm{a-b}_{P_{\Sigma}}^2.
\end{equation}
The following lemma shows that $M$ is positive semidefinite
and that $\mathcal T_\gamma$ is nonexpansive in the induced seminorm. The range-space construction of Section~\ref{sec:preconditioning} then yields a nonexpansive
operator on the Hilbert space associated with $\ran M$.

\begin{lemma}\label{lem:affine-product}
Suppose that Assumption~\ref{ass:affine} and
\eqref{eq:affine-bilinear} hold.
Then the following statements hold:
\begin{enumerate}
\renewcommand{\labelenumi}{(\roman{enumi})}
\item
The matrix $M$ is positive semidefinite, with
$\ker M=\ker P_{\Sigma}\times\ker P_{\Sigma}$.
\item
For every $u,v,\hat u,\hat v\in\R^n$, we have 
\begin{equation}
    \label{ineq:nonexp4Tgamma}
    \norm{\mathcal T_\gamma(u,v)-\mathcal T_\gamma(\hat u,\hat v)}_{M}
\le
\norm{(u,v)-(\hat u,\hat v)}_{M}.
\end{equation}
\item
Set $\Wcal:=\ran M=(\ran P_{\Sigma})^2$ and
$M_{\Wcal}:=M|_{\Wcal}$.
Equip $\Wcal$ with the inner product
\begin{equation*}\label{eq:affine-range-inner-product}
\ip{\xi}{\eta}_{\Wcal}:=\ip{M_{\Wcal}^{-1}\xi}{\eta},
\qquad \xi,\eta\in \Wcal,
\end{equation*}
and let $\norm{\cdot}_{\Wcal}$ denote the induced norm. Then
$\widetilde{\mathcal T}_\gamma\bigl(M(u,v)\bigr)
:=M\mathcal T_\gamma(u,v)$ defines a nonexpansive operator on $(\Wcal,\ip{\cdot}{\cdot}_{\Wcal})$.
Moreover, $s\mapsto M(s,s)$ is a bijection
from $\Scal$ onto $\Fix\widetilde{\mathcal T}_\gamma$.

\end{enumerate}
\end{lemma}

\begin{proof}
\noindent(i)
Since $P+\gamma\Sigma\succeq P_{\Sigma}$, the definition
\eqref{eq:affine-product-metric} and the step-size condition
\eqref{eq:affine-bilinear} give, for all $a,b\in\R^n$,
\[
\begin{aligned}
\ip{(a,b)}{M(a,b)}
&\ge \norm{a}_{P_{\Sigma}}^2-\norm{a}_{P_{\Sigma}}\norm{a-b}_{P_{\Sigma}}
     +\frac12\norm{a-b}_{P_{\Sigma}}^2\\
&=\frac12\norm{a}_{P_{\Sigma}}^2
  +\frac12\bigl(\norm{a}_{P_{\Sigma}}-\norm{a-b}_{P_{\Sigma}}\bigr)^2
\ge0.
\end{aligned}
\]
It follows that $M\succeq0$ and $\ker M\subseteq\ker P_{\Sigma}\times\ker P_{\Sigma}$.
Conversely, if $a,b\in\ker P_{\Sigma}$, then all terms in
\eqref{eq:affine-product-metric} vanish by
\eqref{inclu: KP} and~\eqref{eq: AZU=0}.
Thus, $\ker M=\ker P_{\Sigma}\times\ker P_{\Sigma}$.

\noindent(ii)
Let $(p,u)=\mathcal T_\gamma(u,v)$ and
$(\hat p,\hat u)=\mathcal T_\gamma(\hat u,\hat v)$,
and write
\[
\Delta u:=u-\hat u,\qquad
\Delta v:=v-\hat v,\qquad
\Delta p:=p-\hat p.
\]
The two backward inclusions and the monotonicity
of $F+B-\Sigma$ give
\begin{equation}\label{eq:affine-backward-comparison}
\ip{\Delta p}{P(\Delta u-\Delta p)}
-\gamma\ip{\Delta p}{A(2\Delta u-\Delta v-\Delta p)}
\ge\gamma\norm{\Delta p}_\Sigma^2.
\end{equation}
Expanding~\eqref{eq:affine-product-metric} and using
\eqref{eq:affine-backward-comparison} together with
$P_{\Sigma}-P=\gamma\Sigma/2$, we obtain
\[
\begin{aligned}
\norm{(\Delta p,\Delta u)}_{M}^2
-\norm{(\Delta u,\Delta v)}_{M}^2
&\le
-\frac12\norm{\Delta u-\Delta p}_{P_{\Sigma}}^2
-\frac12\norm{\Delta u-\Delta v}_{P_{\Sigma}}^2\\
&\quad
+2\gamma\ip{\Delta u-\Delta p}{A(\Delta u-\Delta v)}
-\frac{\gamma}{2}\norm{\Delta u+\Delta p}_\Sigma^2\le 0,
\end{aligned}
\]
where the last inequality follows from
\eqref{eq:affine-bilinear} and Young's inequality.
This proves~\eqref{ineq:nonexp4Tgamma}.

\noindent(iii)
By part~(i), the range-space construction of
Section~\ref{sec:preconditioning}, applied to $M$,
makes $\Wcal$ a Hilbert space and gives
\begin{equation}\label{eq:affine-range-isometry}
\norm{M h}_{\Wcal}=\norm{h}_{M},
\qquad h\in\R^{2n}.
\end{equation}
Part~(ii) therefore shows that
$\widetilde{\mathcal T}_\gamma$ is well defined
and nonexpansive.

Let $\xi=M(u,v)$ be a fixed point, and write
$(p,u)=\mathcal T_\gamma(u,v)$.
Then $(p-u,u-v)\in\ker M=(\ker P_{\Sigma})^2$.
Since $P$ and $A$ vanish on $\ker P_{\Sigma}$, we have
$Pu=Pp$ and $F(2u-v)=F(p)$.
The backward inclusion thus gives $p\in \Scal$,
and $\xi=M(p,p)$.
Conversely, every $s\in \Scal$ satisfies
$\mathcal T_\gamma(s,s)=(s,s)$, so
$M(s,s)\in\Fix\widetilde{\mathcal T}_\gamma$.
Finally, if $M(s,s)=M(t,t)$
for $s,t\in \Scal$, then $s-t\in\ker P_{\Sigma}$ and hence
$Ps-\gamma F(s)=Pt-\gamma F(t)$.
Assumption~\ref{ass:affine}(ii) gives $s=t$.

\end{proof}

\subsection{Residual Bounds and Convergence}
We now establish a natural-residual bound and convergence
of the SHARP iterates.

\begin{theorem}\label{thm:affine-convergence}
Suppose that Assumption~\ref{ass:affine}
and the step-size condition~\eqref{eq:affine-bilinear}
hold.
Let $\{(x_k,y_k)\}_{k\ge0}$ be generated by the update
formulas of Algorithm~\ref{alg:darg} from $x_0\in\R^n$.
Then the following statements hold:
\begin{enumerate}
\renewcommand{\labelenumi}{(\roman{enumi})}
\setlength{\itemsep}{0pt}
\setlength{\parsep}{0pt}

\item
The natural residual satisfies
\begin{equation}\label{eq:affine-rate}
\rnat(x_k)
\le
\frac{2\sqrt{\norm{P_{\Sigma}}}}{\gamma k}
\inf_{s\in \Scal}\norm{x_0-s}_{P+\gamma\Sigma},
\qquad k\ge1.
\end{equation}

\item
The sequences $\{x_k\}_{k \ge 0}$ and $\{y_k\}_{k \ge 0 }$ converge to
some $\bar x\in \Scal$ in the Euclidean norm.
Moreover, $\bar x$ is the unique minimizer of
\begin{equation}\label{eq:affine-selected-solution}
\min_{x\in \Scal}\norm{x-x_0}_P^2.
\end{equation}
\end{enumerate}
\end{theorem}

\begin{proof}
Set $\xi_k:=M(u_k,v_k)$ for $k\ge0$.
By Lemmas~\ref{lem:affine-well-posedness}(ii)
and~\ref{lem:affine-product}(iii),
\begin{equation*}\label{eq:affine-range-halpern}
\xi_{k+1}
=\frac{1}{k+2}\xi_0
+\frac{k+1}{k+2}\widetilde{\mathcal T}_\gamma(\xi_k),
\qquad \xi_0=M(x_0,x_0).
\end{equation*}
Moreover,~\eqref{eq:affine-product-metric}
and~\eqref{eq:affine-range-isometry} give
\begin{equation}\label{eq:affine-anchor-distance}
\norm{\xi_0-M(s,s)}_{\Wcal}
=\norm{x_0-s}_{P+\gamma\Sigma},
\qquad s\in \Scal.
\end{equation}

\noindent(i)
Fix $k\ge1$. By~\eqref{eq:affine-output}, we have 
$(x_k,u_{k-1})=\mathcal T_\gamma(u_{k-1},v_{k-1})$.
Lieder's bound~\cite[Theorem~2.1]{lieder2021convergence}, together with
\eqref{eq:affine-anchor-distance}, yields
\[
\begin{aligned}
R&:=\norm{\xi_{k-1}
-\widetilde{\mathcal T}_\gamma(\xi_{k-1})}_{\Wcal}
=\norm{(u_{k-1},v_{k-1})-\mathcal T_\gamma(u_{k-1},v_{k-1})}_{M}\le\frac{2}{k}\inf_{s\in \Scal}
\norm{x_0-s}_{P+\gamma\Sigma}.
\end{aligned}
\]
Set $a:=u_{k-1}-x_k$ and $h:=y_{k-1}-x_k$. Then we have $R=\|(a,a-h)\|_M$ and $d_k=\gamma^{-1}Pa-Ah\in(F+B)(x_k)$.
For any $z\in\R^n$, substituting the expressions for $R^2$ and
$\gamma d_k$ yields
\[
\begin{aligned}
R^2-2\gamma\ip{d_k}{z}+\norm{z}_{P_\Sigma}^2
&=\norm{a-z}_P^2+\gamma\norm{a}_\Sigma^2+\frac{\gamma}{2}\norm{z}_\Sigma^2+\frac12\norm{h}_{P_\Sigma}^2-2\gamma\ip{a-z}{Ah}\\
&\ge\frac12\norm{a-z}_P^2+\frac{\gamma}{2}\norm{a}_\Sigma^2
  +\frac{\gamma}{4}\norm{a+z}_\Sigma^2
\ge0.
\end{aligned}
\]
The first inequality follows from~\eqref{eq:affine-bilinear}, which gives
$2\gamma\ip{a-z}{Ah}\le \|a-z\|_{P_\Sigma}\|h\|_{P_\Sigma}\le
\frac12\bigl(\norm{a-z}_{P_\Sigma}^2+\norm{h}_{P_\Sigma}^2\bigr)$. Taking $z=\gamma d_k/\norm{P_{\Sigma}}$ and using
$\norm{z}_{P_{\Sigma}}^2\le\norm{P_{\Sigma}}\norm{z}^2$ gives
$\norm{d_k}\le\sqrt{\norm{P_{\Sigma}}}R/\gamma$. Combining this estimate with $\rnat(x_k)\le\norm{d_k}$ yields~\eqref{eq:affine-rate}.

\noindent(ii)
By Wittmann's theorem~\cite[Theorem~2]{wittmann1992approximation},
$\xi_k$ converges in $\Wcal$ to the metric projection
$\bar\xi$ of $\xi_0$ onto $\Fix\widetilde{\mathcal T}_\gamma$.
Lemma~\ref{lem:affine-product}(iii) gives
$\bar\xi=M(\bar x,\bar x)$ for a unique $\bar x\in\Scal$.
By~\eqref{eq:affine-anchor-distance}, $\bar x$ uniquely minimizes
$\norm{x_0-s}_{P+\gamma\Sigma}^2$ over $s\in\Scal$.
For $s,t\in\Scal$, monotonicity of $F+B-\Sigma$ gives
$0\le-\norm{s-t}_\Sigma^2$.
Since $\Sigma\succeq0$, we have $\Sigma s=\Sigma t$,
so $\norm{x_0-s}_\Sigma^2$ is constant on $\Scal$.
Thus, $\bar x$ is also the unique minimizer
of~\eqref{eq:affine-selected-solution}.

Since $P$ and $A$ vanish on $\ker P_\Sigma$ and
$\ker M=(\ker P_\Sigma)^2$, the backward input
$Pu-\gamma F(2u-v)$ depends only on $M(u,v)$ and induces a
continuous affine map on $\Wcal$. Applying this map to
$\xi_k\to\bar\xi=M(\bar x,\bar x)$ gives
\[
q_k:=Pu_k-\gamma F(2u_k-v_k)
\longrightarrow
\bar q:=P\bar x-\gamma F(\bar x).
\]

Let $G:=P+\gamma B$. Since $P$ is continuous and monotone
and $B$ is maximal monotone, the sum theorem~\cite{bauschke2017convex}
implies that $G$ and its inverse are maximal monotone.
The backward update gives $x_{k+1}\in G^{-1}(q_k)$.
Since $\bar x\in\Scal$, Lemma~\ref{lem:affine-well-posedness}(i)
gives $G^{-1}(\bar q)=\{\bar x\}$.
By~\cite[Theorem~12.63(c)]{rockafellar1998variational},
$G^{-1}$ is continuous at $\bar q$, and hence $x_{k+1}\to\bar x$.
The reflected update~\eqref{eq:darg-reflected-update}
then gives $y_k\to\bar x$.
\end{proof}

\begin{remark}\label{rem:affine-comparison}
Suppose that $F$ is affine on $\R^n$ and
Assumption~\ref{ass:degenerate} holds. With $\Sigma=0$, $F+B$ is
monotone and $\ran P_\Sigma-\{\gamma c\}=\ran P$ since
$c=F(0)\in\ran P$. Hence Assumption~\ref{ass:affine} holds.
By Lemma~\ref{lem:exact-reduction}(i) and~\eqref{eq:range-metric-identities},
\[
\lvert\ip{Au}{v}\rvert
=\bigl\lvert\ip{F_{\Vcal}(Pu)-F_{\Vcal}(0)}{Pv}_{\Vcal}\bigr\rvert
\le L_P\norm{u}_P\norm{v}_P,
\qquad u,v\in\R^n.
\]
Thus~\eqref{eq:affine-bilinear} holds for $0<\gamma\le1/(2L_P)$.
In this range, Theorem~\ref{thm:affine-convergence} gives the second upper bound in~\eqref{eq:darg-rate} and convergence in the Euclidean norm,
including at $\gamma=1/(2L_P)$.

To illustrate the relaxed range condition, take
\[
P=\operatorname{diag}(1,0),\qquad
A=\operatorname{diag}(0,1),\qquad
B(x)=\operatorname{diag}(0,3)x,\qquad c=0.
\]
With $\Sigma=\operatorname{diag}(0,4)$, we have $F+B-\Sigma=0$,
$P_\Sigma=\operatorname{diag}(1,2\gamma)$, and
$P+\gamma B=\operatorname{diag}(1,3\gamma)$.
Assumption~\ref{ass:affine} and~\eqref{eq:affine-bilinear} hold
for every $\gamma>0$. Since $\ran F=\{0\}\times\R$ and $\ran P=\R\times\{0\}$, we have
$\ran F\not\subseteq\ran P$.
\end{remark}

\section{Structure-Exploiting Preconditioner Design}
\label{sec:preconditioner-design}

The preceding convergence analysis provides a general framework
for SHARP, offering flexibility in the choice of preconditioners.
To translate this flexibility into computational efficiency,
the preconditioner should be tailored to the structure of the
problem. We illustrate this principle through two examples:
a quadratic recourse game using a preconditioner derived
from the sGS decomposition and a multicommodity-flow linear
program using a block preconditioner with Schur-complement
elimination. Both problems are considered in the numerical
experiments of Section~\ref{sec:numerics}.

\subsection{An sGS-Based Preconditioner for Quadratic Recourse Games}
\label{subsec:recourse-sgs}

We consider a two-stage quadratic game with $p$ players, $m$ locations, and $S$ scenarios, each with probability $1/S$, following the framework in~\cite[Section~4]{chen2019twostagesge}. Each player $i\in\{1,\ldots,p\}$ chooses a first-stage production vector $x_i\in[0,\bar x]^m$ before the scenario is observed, and a second-stage adjustment vector $y_{is}\in[-b,b]^m$ once scenario $s\in\{1,\ldots,S\}$ is realized, where $\bar x,b>0$. Given the production plans and scenario adjustments of the other players, player $i$ solves
\begin{equation*}\label{eq:recourse-player}
\begin{aligned}
\min_{\substack{x_i\in[0,\bar x]^m\\
                  y_{is}\in[-b,b]^m,\ s=1,\ldots,S}}\quad
&q_i^\top x_i+\frac1S\sum_{s=1}^S
\left[\frac\omega2\|D_{is}x_i+Ry_{is}-a_{is}\|^2
      +\frac\kappa2\|y_{is}\|^2\right]\\
&+\sum_{j\ne i}(K_{\mathrm{int}})_{ij}\left[
\langle x_i,x_j\rangle+
\frac1S\sum_{s=1}^S\langle y_{is},y_{js}\rangle\right].
\end{aligned}
\end{equation*}
In scenario $s$, $D_{is}x_i$ is the scenario-dependent production and $a_{is}\in\mathbb R^m$ is the demand target, with $D_{is}\succ0$ diagonal. 
The recourse matrix $R\in\mathbb R^{m\times m}$ has unit diagonal and entries $-1$ on its first subdiagonal, yielding $Ry=(y_1,y_2-y_1,\ldots,y_m-y_{m-1})^\top$. For a chain of $m$ locations, $Ry$ gives the net supply adjustment at each location, with $y_1,\ldots,y_{m-1}$ capturing bidirectional transfers between adjacent locations and $y_m$ representing an external adjustment at location $m$.
The vector $q_i\in\mathbb R_{++}^m$ specifies the production costs, and $\omega,\kappa>0$ weight the demand-mismatch and adjustment costs, respectively. The skew-symmetric matrix $K_{\mathrm{int}}\in\mathbb R^{p\times p}$ describes interactions among players.

Let $x$ and $y_s$ denote the stacked production and scenario-adjustment vectors of all players, respectively. Set $n=pm$ and $N=(S+1)n$, and introduce the scaled scenario variables $Y_s=y_s/\sqrt S$ for $s=1,\ldots,S$. This scaling removes the factor $1/S$ from the quadratic and interaction terms involving only the scenario adjustments. Writing $u=(x;Y_1;\ldots;Y_S)\in\mathbb R^N$, we define
\[
D_s=\operatorname{Diag}(D_{1s},\ldots,D_{ps}),\qquad
\mathcal R={\mathrm I}_p\otimes R,\qquad
\mathcal K={\mathrm I}_{S+1}\otimes(K_{\mathrm{int}}\otimes {\mathrm I}_m).
\]
Since each player's objective is convex in its own decision variables, the joint first-order optimality conditions characterize a Nash equilibrium and take the form of the monotone inclusion
\begin{equation}\label{eq:recourse-design-system}
0\in G(u)+\mathcal N_{\mathcal X}(u),\qquad
G(u)=\mathcal Ku+Qu+c,
\end{equation}
where $\mathcal X=[0,\bar x]^n\times[-b/\sqrt S,b/\sqrt S]^{Sn}$. With $Y=(Y_1;\ldots;Y_S)$, the quadratic and affine terms are
\[
Q:=\begin{pmatrix}Q_{xx}&E\\E^\top&Q_{YY}\end{pmatrix},\qquad
c:=(c_x;c_{Y_1};\ldots;c_{Y_S}),\qquad
E:=\frac{\omega}{\sqrt{S}}(D_1^\top\mathcal R,\ldots,D_S^\top\mathcal R),
\]
with
\[
\begin{aligned}
Q_{xx}&:=\frac{\omega}{S}\sum_{s=1}^S D_s^\top D_s, &\qquad
Q_{YY}&:={\mathrm I}_S\otimes(\kappa{\mathrm I}_n+\omega\mathcal R^\top\mathcal R),\\
c_x&:=q-\frac{\omega}{S}\sum_{s=1}^S D_s^\top a_s, &\qquad
c_{Y_s}&:=-\frac{\omega}{\sqrt{S}}\mathcal R^\top a_s,
\end{aligned}
\]
where $q$ and $a_s$ stack the corresponding player vectors. Combining $\kappa>0$ and $D_s\succ0$ with the identity $u^\top Qu=\kappa\|Y\|^2+\omega\sum_{s=1}^S\|D_sx/\sqrt S+\mathcal RY_s\|^2$ shows that $Q\succ0$. Since $\mathcal K$ is skew-symmetric, $G$ is continuous and strongly monotone. Because $\mathcal X$ is nonempty, compact, and convex, the inclusion~\eqref{eq:recourse-design-system} admits a unique solution $u_\star$, corresponding to the unique Nash equilibrium of the game.

A direct resolvent evaluation of $Qu+c+\mathcal N_{\mathcal X}(u)$ in~\eqref{eq:recourse-design-system} requires solving a box-constrained quadratic problem. Although it separates across players, the block $E$ couples each player's production plan $x_i$ with its adjustments across all scenarios. Even with $x_i$ fixed, each scenario subproblem retains box constraints and a tridiagonal Hessian ${\mathrm I}_m+\gamma(\kappa{\mathrm I}_m+\omega R^\top R)$. Following the variable-splitting technique in~\cite[Section~3]{sun2026sgs}, we separate the box constraint from the quadratic coupling by introducing an auxiliary variable $z\in\mathcal Z:=\sqrt\beta\mathcal X$ with $\beta>0$ and enforcing $z=\sqrt\beta u$ via a multiplier $\lambda\in\mathbb R^N$. Setting $w=(z;u)$ and $v=(w;\lambda)$, we obtain the lifted inclusion $0\in F(v)+B(v)$, where
\begin{equation}\label{eq:recourse-design-splitting}
F(v)=\begin{pmatrix}0\\\mathcal Ku\\0\end{pmatrix},
\qquad
B(v)=\begin{pmatrix}
\mathcal N_{\mathcal Z}(z)-\lambda\\
Qu+c+\sqrt\beta\lambda\\
z-\sqrt\beta u
\end{pmatrix}.
\end{equation}
Eliminating $z$ and $\lambda$ recovers~\eqref{eq:recourse-design-system}, so the two inclusions are equivalent. The lifted inclusion~\eqref{eq:recourse-design-splitting} has the unique solution
\[
v_\star=(\sqrt\beta u_\star;\,u_\star;\,-G(u_\star)/\sqrt\beta).
\]
Here, $F$ is linear and skew-symmetric, while $B$ is the sum of the subdifferential of the proper closed convex function $\delta_{\mathcal Z}(z)+\tfrac12u^\top Qu+c^\top u$ and a continuous linear skew-symmetric operator, so $B$ is maximal monotone.

Fix $\beta,\tau,\gamma>0$ and set
$C=(-{\mathrm I}_N,\sqrt{\beta}{\mathrm I}_N)$ and
$\bar Q=\operatorname{Diag}(0,Q)$. Define
\[
H:=\operatorname{Diag}(0,\tau{\mathrm I}_N)+\gamma(\bar Q+C^\top C).
\]
The quadratic form
$w^\top Hw=\tau\|u\|^2+\gamma u^\top Qu+
\gamma\|\sqrt{\beta}u-z\|^2$ shows that $H\succ0$.
Decompose $H=D+U+U^\top$ in the three-block order $(z,x,Y)$, where $D$ is block diagonal and $U$ is strictly block upper triangular:
\begin{equation*}\label{eq:recourse-sgs-diagonal-blocks}
D=\operatorname{Diag}(\gamma{\mathrm I}_N,D_x,D_Y),\qquad
D_x=\tau{\mathrm I}_n+\gamma(Q_{xx}+\beta{\mathrm I}_n),\qquad
D_Y=\tau{\mathrm I}_{Sn}+\gamma(Q_{YY}+\beta{\mathrm I}_{Sn}).
\end{equation*}
Following the sGS decomposition in~\cite{li2019sgs}, we introduce the sGS correction $T:=UD^{-1}U^\top$ and the regularized matrix $H_{\mathrm{sGS}}:=H+T$, and define the preconditioner by
\begin{equation}\label{eq:recourse-sgs-metric}
P_{\mathrm{sGS}}=
\begin{pmatrix}
H_{\mathrm{sGS}}-\gamma\bar Q&-\gamma C^\top\\
-\gamma C&\gamma{\mathrm I}_N
\end{pmatrix}.
\end{equation}

Let $r=(r_z;r_u;r_\lambda)$ and $v^+=(w^+;\lambda^+)$. In the
block of $r\in(P_{\mathrm{sGS}}+\gamma B)(v^+)$ corresponding
to $w^+$, the terms $-\gamma\bar Qw^+$ and
$-\gamma C^\top\lambda^+$ in $P_{\mathrm{sGS}}v^+$ cancel
the corresponding terms in $\gamma B(v^+)$. Setting $h_z=r_z$
and $(h_x;h_Y)=r_u-\gamma c$, we obtain the equivalent system
\begin{equation}\label{eq:recourse-sgs-backward-system}
(h_z;h_x;h_Y)\in H_{\mathrm{sGS}}w^+
+\gamma\mathcal N_{\mathcal Z\times\mathbb R^N}(w^+),\qquad\lambda^+=r_\lambda/\gamma+2Cw^+.
\end{equation}
Since $H_{\mathrm{sGS}}=H+UD^{-1}U^\top=(D+U)D^{-1}(D+U^\top)\succ0$ and the normal cone acts only on the first block $z$, the sGS decomposition theorem~\cite[Theorem~1]{li2019sgs} allows $w^+$ to be computed exactly by a backward sweep over $(Y,x)$, a box projection onto $\mathcal Z$, and a forward sweep over $(x,Y)$. Partitioning $z^+=(z_x^+;z_Y^+)$ conformably with $u^+=(x^+;Y^+)$, we obtain
\begin{equation*}\label{eq:recourse-sgs-sweep}
\begin{aligned}
\widetilde Y&=D_Y^{-1}h_Y,\qquad
\widetilde x=D_x^{-1}(h_x-\gamma E\widetilde Y),\qquad
z^+=\Pi_{\mathcal Z}\bigl(h_z/\gamma+\sqrt\beta(\widetilde x;\widetilde Y)\bigr),\\
x^+&=\widetilde x+\gamma\sqrt\beta\,D_x^{-1}z_x^+,\qquad
Y^+=\widetilde Y+\gamma D_Y^{-1}
\bigl(\sqrt\beta\,z_Y^+-E^\top x^+\bigr). 
\end{aligned}
\end{equation*}
The multiplier $\lambda^+$ is then recovered from~\eqref{eq:recourse-sgs-backward-system}. Note that $D_x$ is diagonal, while $D_Y$ consists of $pS$ identical positive-definite tridiagonal blocks
$(\tau+\gamma(\beta+\kappa)){\mathrm I}_m+\gamma\omega R^\top R$. 
Each generalized backward step therefore requires one box projection, two diagonal solves with $D_x$, two batches of tridiagonal solves with $D_Y$, and sparse matrix--vector products. The common tridiagonal factorization is reused across players, scenarios, and iterations.

We now verify the assumptions of Theorem~\ref{thm:affine-convergence} with $P=P_{\mathrm{sGS}}$ and $\Sigma=0$. For $v=(z;u;\lambda)$ with $z=(z_x;z_Y)$, a direct expansion of~\eqref{eq:recourse-sgs-metric} yields
\begin{equation}\label{eq:recourse-metric-energy}
\begin{aligned}
\|v\|_{P_{\mathrm{sGS}}}^2
&=\tau\|u\|^2+\gamma\|\lambda-Cw\|^2
  +\gamma^2\beta\|z_x\|_{D_x^{-1}}^2
  +\gamma^2\|\sqrt\beta z_Y-E^\top x\|_{D_Y^{-1}}^2.
\end{aligned}
\end{equation}
If this expression vanishes, the first term gives $u=0$, the last two terms give $z=0$, and the second term then gives $\lambda=0$, which implies $P_{\mathrm{sGS}}\succ0$. Since $B$ is maximal monotone,
$(P_{\mathrm{sGS}}+\gamma B)^{-1}$ is everywhere defined and
single-valued. Together with the monotonicity of $F+B$ and the
existence of $v_\star$, this verifies Assumption~\ref{ass:affine}
with $\Sigma=0$. Moreover, using $\|v\|_{P_{\mathrm{sGS}}}^2\ge\tau\|u\|^2$ from~\eqref{eq:recourse-metric-energy} and $\|\mathcal K\|_2=\|K_{\mathrm{int}}\|_2$, we obtain for all $v=(z;u;\lambda)$ and $\widehat v=(\widehat z;\widehat u;\widehat\lambda)$,
\begin{equation*}\label{eq:recourse-forward-bound}
|\langle F(v),\widehat v\rangle|
=|\langle\mathcal Ku,\widehat u\rangle|
\le\|\mathcal K\|_2\|u\|\|\widehat u\|
\le L_P\|v\|_{P_{\mathrm{sGS}}}\|\widehat v\|_{P_{\mathrm{sGS}}},
\qquad
L_P:=\frac{\|K_{\mathrm{int}}\|_2}{\tau}.
\end{equation*}
Thus, Theorem~\ref{thm:affine-convergence} applies whenever
$ 2\gamma\|K_{\mathrm{int}}\|_2\le\tau $.

\subsection{A Block Preconditioner with Schur-Complement Elimination}
\label{subsec:network-schur}

As a second example, we consider a multicommodity-flow problem on a directed graph with $N$ nodes and $E$ edges, whose underlying undirected graph is connected. Let $D\in\R^{N\times E}$ be the incidence matrix. Since the graph is connected, $\ker D^{\top}=\operatorname{span}\{\mathbf1_N\}$ and $\ran D=\mathbf1_N^\perp$. Thus, the graph Laplacian $\mathrm L_G:=DD^{\top}$ is singular on $\R^N$ but positive definite on $\mathbf1_N^\perp$. For commodity $i=1,\ldots,K$, let $f_i\in\R^E$, $d_i\in\mathbf1_N^\perp$, and $c_i\in\R^E$ denote its flow, demand, and cost vectors, respectively. With a common edge-capacity vector $a>0$ and a nonnegative slack vector $s$ representing unused capacity, the problem is
\begin{equation}\label{eq:flow-lp}
\min\left\{
\sum_{i=1}^K\ip{c_i}{f_i}:
Df_i=d_i\ (i=1,\ldots,K),\quad
\sum_{i=1}^K f_i+s=a,\quad
f_i,s\ge0
\right\}.
\end{equation}

We assume that the LP~\eqref{eq:flow-lp} is feasible. The capacity constraints and nonnegativity imply $0\le f_i,s\le a$ componentwise; hence imposing these box bounds does not change the feasible set. Set
\[
x=(f_1;\ldots;f_K;s),\quad
\ell=0,\quad
u=(a;\ldots;a),\quad
c=(c_1;\ldots;c_K;0),\quad
b=(d_1;\ldots;d_K;a),
\]
where $u$ contains $K+1$ copies of $a$, and define
\[
A=
\begin{pmatrix}
D & & & 0\\
  & \ddots & & \vdots\\
  & & D & 0\\
{\mathrm I}_E & \cdots & {\mathrm I}_E & {\mathrm I}_E
\end{pmatrix}.
\]
Then, the LP~\eqref{eq:flow-lp} is equivalently written as
\[
\min_{x\in[\ell,u]}\ip{c}{x}
\qquad\text{subject to}\qquad Ax=b.
\]

Let $y=(y_1;\ldots;y_K;\lambda)$ denote the multipliers associated with the flow-balance and capacity constraints. Since $D^\top\mathbf1_N=0$, adding a constant multiple of $\mathbf1_N$ to any $y_i$ leaves the KKT conditions unchanged. We remove this ambiguity by restricting the multipliers to $\mathcal Y:=(\mathbf1_N^\perp)^K\times\R^E$.
Indeed, $\ran D=\mathbf1_N^\perp$ and the block structure of $A$ imply $\ran A=\mathcal Y$ and $\ker A^\top=\mathcal Y^\perp$. The KKT conditions can therefore be written in the normalized form
\begin{equation*}\label{inclusion: prob}
0\in F(x,y)+B(x,y),\qquad
F(x,y)=
\begin{pmatrix}-A^\top y\\Ax\end{pmatrix},\qquad
B(x,y)=
\begin{pmatrix}
c+\mathcal N_{[\ell,u]}(x)\\
-b+\mathcal N_{\mathcal Y}(y)
\end{pmatrix}.
\end{equation*}
Here, $F$ is linear and skew-symmetric, and $B$ is maximally monotone. Since $Ax-b\in\mathcal Y$ and $\mathcal N_{\mathcal Y}(y)=\mathcal Y^\perp$ for $y\in\mathcal Y$, the second inclusion enforces precisely $Ax=b$ and $y\in\mathcal Y$, preserving the original KKT conditions up to the multiplier normalization.

We use a scalar primal block to obtain a box projection in the
backward step and choose the dual block to measure the multiplier
through the coupling term $A^\top y$. This gives
\begin{equation*}\label{eq:flow-schur-metric}
P_\beta:=\operatorname{Diag}
\bigl(\beta{\mathrm I}_{(K+1)E},\beta^{-1}AA^\top\bigr),
\qquad \beta>0.
\end{equation*}
The induced seminorm satisfies $\norm{(x,y)}_{P_\beta}^2
=\beta\norm{x}^2+\beta^{-1}\norm{A^\top y}^2$.
Since $\ker(AA^\top)=\mathcal Y^\perp$, the preconditioner
$P_\beta$ is only positive semidefinite.

For any $\gamma>0$ and input $r=(r_x,r_y)$, the generalized backward step $(x^+,y^+)=(P_\beta+\gamma B)^{-1}r$ takes the form
\begin{equation}\label{eq:flow-generalized-resolvent}
x^+=\Pi_{[\ell,u]}
\left(\frac{r_x-\gamma c}{\beta}\right),\text{ and }AA^{\top}y^+=\beta\bigl(\Pi_{\mathcal Y}r_y+\gamma b\bigr),\quad y^+\in\mathcal Y.
\end{equation}
Since $AA^\top$ is positive definite on $\mathcal Y$, the backward step is uniquely defined for every input despite the degeneracy of $P_\beta$. Moreover, the dual component of each backward input in Algorithm~\ref{alg:darg} belongs to $\ran A=\mathcal Y$, so the projection $\Pi_{\mathcal Y}$ can be omitted.

We next verify the convergence conditions. Since the feasible set of~\eqref{eq:flow-lp} is nonempty and compact, linear programming duality guarantees the existence of a KKT point. Together with the maximal monotonicity of $B$ and the well-defined backward step~\eqref{eq:flow-generalized-resolvent}, this verifies Assumption~\ref{ass:affine}. Moreover, for $w=(x;y)$ and $\widehat w=(\widehat x;\widehat y)$, the weighted Cauchy--Schwarz inequality gives
\[
\bigl|\ip{F(w)}{\widehat w}\bigr|
=\bigl|\ip{x}{A^\top\widehat y}
-\ip{A^\top y}{\widehat x}\bigr|
\le\norm{w}_{P_\beta}\norm{\widehat w}_{P_\beta}.
\]
Thus, the step-size condition~\eqref{eq:affine-bilinear} holds for $0<\gamma\le1/2$, independently of $\beta$. For fixed $\beta>0$ and $0<\gamma\le1/2$, Theorem~\ref{thm:affine-convergence} establishes an $\mathcal O(1/k)$ last-iterate natural-residual bound for the normalized KKT inclusion and strong convergence of the unrestarted SHARP iterates to a normalized KKT point, including at the boundary $\gamma=1/2$.

We now derive an efficient Schur-complement realization of the generalized backward step~\eqref{eq:flow-generalized-resolvent}. Given an input $r=(r_x;r_y)$, write $r_y=(r_1;\ldots;r_K;r_\lambda)$ and set
\[
g_i=\beta\bigl(\Pi_{\mathbf1_N^\perp}r_i+\gamma d_i\bigr),
\quad i=1,\ldots,K,
\qquad
t=\beta(r_\lambda+\gamma a).
\]
The dual system in~\eqref{eq:flow-generalized-resolvent} then takes the block form
\[
\mathrm L_Gy_i^++D\lambda^+=g_i,
\quad i=1,\ldots,K,
\qquad
D^\top\sum_{i=1}^K y_i^++(K+1)\lambda^+=t.
\]
Applying $D$ to the second equation and subtracting the sum of the first $K$ equations gives
$D\lambda^+=Dt-\sum_{i=1}^K g_i$.
Substitution yields the Schur-complement realization
\begin{equation}\label{eq:flow-arrow-solve}
\mathrm L_Gy_i^+
=g_i+\sum_{j=1}^K g_j-Dt,
\qquad \mathbf1_N^\top y_i^+=0,
\quad i=1,\ldots,K\text{  and  }\lambda^+=\frac{1}{K+1}
\left(t-D^\top\sum_{i=1}^K y_i^+\right).
\end{equation}

For each $i$, the right-hand side of the Laplacian equation
in~\eqref{eq:flow-arrow-solve} lies in $\mathbf1_N^\perp$, on which
$\mathrm L_G$ is positive definite. We obtain the unique zero-mean solution by fixing the potential
at a reference node to zero, solving the grounded Laplacian
system, and subtracting the mean of the resulting $N$-vector.
Direct evaluation of~\eqref{eq:flow-arrow-solve} uses $K$
Laplacian solves with a shared factorization.

\makeatletter\setlength{\@fptop}{0pt}\setlength{\@fpsep}{20pt}\makeatother
\newcommand{\tabletime}[3]{#1#2#3}
\newcommand{\tableerror}[2]{\makebox[5.5em][l]{$#1\!\times\!10^{#2}$}}
\newcommand{\tablestatus}[1]{#1}
\makeatletter
\newcommand{\tablepath}[3]{\raisebox{-\dimexpr(\ht\@arstrutbox+\dp\@arstrutbox)*#1/2+#2/2\relax}[0pt][0pt]{#3}}
\makeatother
\section{Numerical Experiments}\label{sec:numerics}

In this section, we evaluate SHARP on three problem families: zero-sum network games, quadratic recourse games, and multicommodity-flow linear programs. For the first two problem families, we compare SHARP with the forward-reflected-backward (FRB) method~\cite{malitsky2020forward} to examine the effects of acceleration and structure-exploiting
preconditioning. For the multicommodity-flow LPs, we implement SHARP using the block preconditioner developed in Section~\ref{subsec:network-schur} and compare the resulting implementation with leading LP solvers HPR-LP-C~\cite{chen2026hprlp} and MOSEK~\cite{mosek2026} to show the computational efficiency
of our method.

We conducted all experiments on a Linux server with Intel Xeon Platinum
8470Q processors, using an allocation of 64 CPU cores and 128~GB of host
memory. SHARP, FRB, and HPR-LP-C ran on an NVIDIA GeForce RTX~5090 GPU
with 32~GB of device memory, whereas MOSEK ran on the CPU. We implemented
SHARP and FRB in C++17 and CUDA~12.8.

For the following experiments, we use SHARP as given in
Algorithm~\ref{alg:darg} with $\gamma L_P=1/2$, where $L_P$ is the
Lipschitz constant of $F_{\Vcal}$ in the range-space norm $\|\cdot\|_{\Vcal}$.
The forward operator $F$ is affine in all three problem families, and
Theorem~\ref{thm:affine-convergence} permits this endpoint for the
preconditioners used below. For the first two families, the FRB
iteration is
\[
x_{k+1}
=(P+\gamma B)^{-1}\bigl(Px_k-\gamma(2F(x_k)-F(x_{k-1}))\bigr),
\qquad x_{-1}=x_0.
\]
We use $\gamma L_P=0.49$ for FRB to satisfy the strict step-size
condition $\gamma L_P<1/2$ in its convergence guarantee
\cite[Corollary~2.6]{malitsky2020forward}, applied in the range space.

To improve the practical performance of SHARP, we use an adaptive
restart strategy based on the norm of the residual vector $d_j$
defined in~\eqref{eq:darg-dk}. Let $j$ denote the iteration index
within the current restart cycle and $k$ the total iteration count.
We use $\|d_1\|$ as the reference value for each cycle and check
for a restart every $N_{\mathrm{check}}=150$ iterations. Following the approach of
Applegate et al.~\cite[Section~3.2]{applegate2021practical}, we use
the following three restart conditions.
\begin{enumerate}
\item[(i)] \emph{Sufficient decrease}:
$\|d_j\|\le\alpha_1\|d_1\|$.
\item[(ii)] \emph{Necessary decay and no local progress}:
$\|d_j\|\le\alpha_2\|d_1\|$ and $\|d_j\|<\|d_{j+N_{\mathrm{check}}}\|$.
\item[(iii)] \emph{Long inner loop}:
$j\ge\alpha_3 k$.
\end{enumerate}
We set $(\alpha_1,\alpha_2,\alpha_3)=(0.1,0.6,0.2)$.
In cases~(i) and~(iii), we restart from the current iterate $x_j$. Condition~(ii) is evaluated after computing $d_{j+N_{\mathrm{check}}}$; if it holds, we restart from the stored iterate $x_j$. In each case, the selected iterate becomes the new anchor.

We assess solution accuracy using relative natural residuals for the two game families and relative feasibility, stationarity, and gap residuals for the LPs. The residuals used in each experiment are defined in the corresponding subsection. The tolerance is set to $10^{-6}$ for all methods. In all three experiments, we impose a time limit of $3600$ seconds per run and for SHARP, FRB and HPR-LP-C we check the stopping criterion every $150$ iterations. 

\subsection{Zero-Sum Network Games}

We first consider zero-sum polymatrix games~\cite{cai2016zerosum} on a star graph with $p$ players. The central player $0$ interacts with players $1,\ldots,p-1$,
and each player chooses a mixed strategy from the simplex
$\Delta_m:=\{v\in\R^m:v\ge0,\ \mathbf1^\top v=1\}$.
For each edge $(0,a)$, let $R_a\in\R^{m\times m}$
be the payoff matrix.
Each player maximizes its own payoff while holding the
other players' strategies fixed.
The central and peripheral players therefore solve,
respectively,
\[
\max_{x_0\in\Delta_m}
\left(-\sum_{b=1}^{p-1}x_0^\top R_bx_b\right),
\qquad
\max_{x_a\in\Delta_m}x_0^\top R_ax_a,
\quad a=1,\ldots,p-1.
\]

Let $x=(x_0,\ldots,x_{p-1})$ and $\mathcal X=\Delta_m^p$.
Since each player's optimization problem is linear over
a simplex, its first-order optimality condition is both
necessary and sufficient.
Combining these conditions shows that finding a Nash
equilibrium is equivalent to solving the monotone inclusion
\begin{equation*}\label{eq:current-game}
0\in F(x)+B(x),
\qquad
F(x)=Gx,
\qquad
B=\mathcal N_{\mathcal X},
\end{equation*}
where
\[
G=
\begin{pmatrix}
0&R_1&\cdots&R_{p-1}\\
-R_1^\top&0&\cdots&0\\
\vdots&\vdots&\ddots&\vdots\\
-R_{p-1}^\top&0&\cdots&0
\end{pmatrix}.
\]
Since $G^\top=-G$, the forward operator $F$ is linear
and monotone, while $B$ is maximal monotone.

The central player participates in $p-1$ interactions,
while each peripheral player participates in only one.
This imbalance motivates a preconditioner based on node degrees.
We compare SHARP and FRB with the identity preconditioner
and the node-degree preconditioner
\[
P_{\mathrm{id}}={\mathrm I}_{pm},
\qquad
P_{\mathrm{deg}}
=\operatorname{Diag}\bigl((p-1){\mathrm I}_m,{\mathrm I}_m,\ldots,{\mathrm I}_m\bigr).
\]
Both are positive definite.
For a block-diagonal preconditioner
$P=\operatorname{Diag}(d_0{\mathrm I}_m,\ldots,d_{p-1}{\mathrm I}_m)$ with $d_a>0$,
the generalized backward step reduces to independent
simplex projections:
\[
\bigl[(P+\gamma\mathcal N_{\mathcal X})^{-1}
       (Pv-\gamma f)\bigr]_a
=
\Pi_{\Delta_m}
\left(v_a-\frac{\gamma}{d_a}f_a\right),
\qquad a=0,\ldots,p-1.
\]
Thus, the two preconditioners require the same projection operations, but assign different effective step-sizes to the players. The identity preconditioner uses $d_a=1$ for every player, whereas the node-degree preconditioner uses $d_0=p-1$ and $d_a=1$ for $a\ge1$.

For either preconditioner, the Lipschitz constant entering the step-size rule is $L_P=\|P^{-\frac{1}{2}}GP^{-\frac{1}{2}}\|_2$.
The block structure of $G$ yields
\[
L_{\mathrm{id}}
=\|G\|_2
=
\sqrt{\lambda_{\max}
\left(\sum_{a=1}^{p-1}R_aR_a^\top\right)},
\qquad
L_{\mathrm{deg}}
=\frac{L_{\mathrm{id}}}{\sqrt{p-1}}.
\]
Both constants are therefore obtained from the largest eigenvalue of the same $m\times m$ Gram matrix. At the same normalized step-size $\gamma L_P$, node-degree preconditioning decreases the coefficient $\gamma/d_a$ of the forward term for the central player by a factor of $\sqrt{p-1}$ and increases those for the peripheral players by the same factor, relative to the identity preconditioner.

For each edge, we generate a matrix with independent standard
Gaussian entries and divide it by a power-method estimate of
its spectral norm to obtain $R_a$.
The problem sizes are listed in Table~\ref{tab:current-zero-sum}.
All methods start from the uniform strategy
$x_{a,0}=\mathbf1/m$ for $a=0,\ldots,p-1$. The relative natural residual used in the stopping criterion is $
\|x-\Pi_{\mathcal X}(x-Gx)\|_2/(
1+\|x\|_\infty+\|Gx\|_\infty)$.

\begin{table}[!htbp]
\centering
\small
\setlength{\tabcolsep}{2.8pt}
\renewcommand{\arraystretch}{1.16}
\caption{Zero-sum network games: running times after warm-up and iteration
counts of SHARP and FRB with identity and node-degree
preconditioners.
T denotes termination at the $3600$-second time limit.
The reported times include the spectral computation for the
chosen preconditioner and one full solver call.}
\label{tab:current-zero-sum}
\vspace{4pt}
\begin{tabular*}{\textwidth}
{@{\extracolsep{\fill}}rrrrrrrrr@{}}
\toprule
Problem size
& \multicolumn{2}{c}{FRB--Identity}
& \multicolumn{2}{c}{SHARP--Identity}
& \multicolumn{2}{c}{FRB--Degree}
& \multicolumn{2}{c}{SHARP--Degree}\\
\cmidrule(lr){1-1}
\cmidrule(lr){2-3}
\cmidrule(lr){4-5}
\cmidrule(lr){6-7}
\cmidrule(lr){8-9}
$(p,m)$
& Time (s) & Iterations
& Time (s) & Iterations
& Time (s) & Iterations
& Time (s) & Iterations\\
\midrule
$(1024,32)$
& 90.80 & 1\,131\,150
& 136.42 & 2\,134\,800
& 12.64 & 158\,700
& \textbf{5.59} & 86\,700\\
$(2048,32)$
& 140.63 & 954\,450
& 246.60 & 1\,861\,800
& 22.09 & 149\,700
& \textbf{7.52} & 56\,700\\
$(4096,32)$
& 385.62 & 1\,385\,250
& 689.04 & 2\,674\,950
& 45.88 & 162\,000
& \textbf{13.85} & 53\,700\\
$(8192,32)$
& 723.48 & 1\,335\,750
& 1238.60 & 2\,612\,550
& 111.40 & 205\,050
& \textbf{17.99} & 37\,800\\
$(16384,32)$
& 2320.80 & 1\,438\,350
& ${\mathrm T}$ & ---
& 299.47 & 185\,100
& \textbf{52.53} & 36\,300\\
$(32768,32)$
& ${\mathrm T}$ & ---
& ${\mathrm T}$ & ---
& 605.83 & 184\,350
& \textbf{120.71} & 40\,650\\
$(65536,32)$
& ${\mathrm T}$ & ---
& ${\mathrm T}$ & ---
& 1266.39 & 194\,550
& \textbf{218.20} & 36\,750\\
\bottomrule
\end{tabular*}
\end{table}

Table~\ref{tab:current-zero-sum} shows that node-degree preconditioning
improves both methods, with particularly large gains for SHARP.
At $(p,m)=(8192,32)$, it reduces the running times of SHARP and FRB
by factors of $68.9$ and $6.5$, respectively. With the node-degree
preconditioner, both methods solve all seven instances, and SHARP
requires fewer iterations and achieves time speedups of
$2.26$--$6.19$ over FRB. The larger gains for SHARP suggest that its
accelerated iteration may be more sensitive to scaling than the FRB
iteration. These results highlight the benefit of combining SHARP with structure-exploiting preconditioning on these
unbalanced networks.

\subsection{Quadratic Recourse Games}
\label{sec:scenario-recourse-sgs}
We next consider the quadratic recourse game introduced in Section~\ref{subsec:recourse-sgs}. We apply SHARP and FRB to the lifted inclusion
\eqref{eq:recourse-design-splitting} and compare block-Jacobi
and sGS preconditioning as the numbers of locations and scenarios increase.

We fix the number of players at $p=128$ and define their interactions through the skew-symmetric matrix
\[
K_{\mathrm{int}}=\begin{pmatrix}
0&W_r\\
-W_r^\top&0
\end{pmatrix},
\qquad r=p/2,
\]
where $W_r\in\mathbb R^{r\times r}$ has unit diagonal and entries $-1$ on its first subdiagonal, with all other entries zero. Consequently, $\|\mathcal K\|_2=\|K_{\mathrm{int}}\|_2=2\cos(\pi/(p+1))$. The location-coupling matrix $R$ retains its unit-difference structure as $m$ increases, and thus increasing $m$ extends the location network rather than refining a fixed spatial domain.

We generate the production efficiencies and demand targets by combining player--location-specific characteristics with scenario-dependent perturbations:
\[
\begin{aligned}
(D_{is})_{jj}&=1+\delta_{ij}+\varepsilon_{is,j},
&\delta_{ij},\varepsilon_{is,j}&\sim\mathcal U[-1/4,1/4],\\
a_{is,j}&=\bar a_{ij}+0.1\xi_{is,j},
&\bar a_{ij}&\sim\mathcal U[0.5,0.9],\quad
\xi_{is,j}\sim\mathcal N(0,1).
\end{aligned}
\]
Here, $\delta_{ij}$ and $\bar a_{ij}$ are sampled once for each player--location pair and shared across scenarios, whereas $\varepsilon_{is,j}$ and $\xi_{is,j}$ are generated independently across scenarios. The production-cost coefficients are given by $q_{i,j}=0.1(0.5+U_{i,j})$, where $U_{i,j}\sim\mathcal U[0,1]$. We fix $\kappa=10^{-5}$, $\omega=10$, $\bar x=1.5$, and $b=1$ across all problem sizes. 
We consider five problem sizes with $m=S=32,64,128,256,512$, increasing the numbers of locations and scenarios simultaneously. The dimension of the original game is $N=pm(S+1)$, ranging from $135\,168$ to $33\,619\,968$. For each size, all four configurations use the same generated instance and are initialized at zero.

For both preconditioners, we set $\beta=0.1$ and
$\tau=2\|\mathcal K\|_2$. Motivated by \cite{borgens2018distributed}, a block-Jacobi preconditioner for solving $0\in F(v)+B(v)$ is obtained by replacing
$H_{\mathrm{sGS}}$ in~\eqref{eq:recourse-sgs-metric} with
\[
H_J=\operatorname{Diag}\left(
2\gamma{\mathrm I}_N,\;
D_x+\gamma(\beta{\mathrm I}_n+Q_{xx}),\;
D_Y+\gamma(\beta{\mathrm I}_{Sn}+Q_{YY}-\kappa{\mathrm I}_{Sn})
\right).
\]
This preconditioner is positive definite and admits
the same Lipschitz bound $L_P=\|\mathcal K\|_2/\tau=1/2$ as sGS. Following the common step-size rule, we use $\gamma=1$ for SHARP and $\gamma=0.98$ for FRB, constructing each preconditioner with the corresponding step-size.
The backward step under this Jacobi metric can be evaluated using one box projection and independent solves of diagonal and tridiagonal systems.

For all four configurations, we evaluate the relative natural residual
$\|u-\Pi_{\mathcal X}(u-G(u))\|_2/(1+\|u\|_\infty+\|G(u)\|_\infty)$ of the unlifted inclusion
\eqref{eq:recourse-design-system} in the scaled variables $u=(x;Y_1;\ldots;Y_S)$.
Each configuration is run in a separate process.
After warm-up, we record the time of a single complete solver call. This time includes solver setup, iterations, and solution recovery, but excludes instance generation and warm-up.

\begin{table}[H]
\centering
\small
\setlength{\tabcolsep}{2.8pt}
\renewcommand{\arraystretch}{1.16}
\caption{Quadratic recourse games: full-call times (s) after warm-up and iteration counts.}
\label{tab:current-recourse}
\vspace{4pt}
\begin{tabular*}{\textwidth}{@{\extracolsep{\fill}}rrrrrrrrr@{}}
\toprule
Problem size
& \multicolumn{2}{c}{FRB--Jacobi}
& \multicolumn{2}{c}{SHARP--Jacobi}
& \multicolumn{2}{c}{FRB--sGS}
& \multicolumn{2}{c}{SHARP--sGS}\\
\cmidrule(lr){1-1}\cmidrule(lr){2-3}\cmidrule(lr){4-5}
\cmidrule(lr){6-7}\cmidrule(lr){8-9}
$(p,m,S)$ & Time (s) & Iterations & Time (s) & Iterations
          & Time (s) & Iterations & Time (s) & Iterations\\
\midrule
$(128,32,32)$   & 2.76    & 107\,100 & 1.37   & 45\,750 & 4.14    & 97\,200  & \textbf{0.69}   & 14\,400\\
$(128,64,64)$   & 16.40   & 248\,400 & 3.26   & 43\,800 & 20.80   & 198\,900 & \textbf{1.48}   & 13\,050\\
$(128,128,128)$ & 113.39  & 347\,400 & 16.30  & 44\,400 & 133.71  & 261\,600 & \textbf{7.49}   & 13\,500\\
$(128,256,256)$ & 504.57  & 340\,800 & 72.84  & 45\,150 & 604.29  & 267\,000 & \textbf{33.55}  & 13\,950\\
$(128,512,512)$ & 2355.91 & 375\,150 & 316.86 & 46\,800 & 2762.19 & 288\,150 & \textbf{149.14} & 14\,700\\
\bottomrule
\end{tabular*}
\end{table}

Table~\ref{tab:current-recourse} shows that SHARP--sGS achieves the
lowest full-call time at every tested size, with speedups of
$1.99$--$2.20$ over SHARP--Jacobi. Its iteration counts remain between
$13\,050$ and $14\,700$ across all five instances. For SHARP, sGS
preconditioning reduces the iteration count by more than a factor of
three relative to block Jacobi, offsetting the cost of the additional
round of $x$- and $Y$-block solves. The FRB results provide a useful
contrast to these gains. Although FRB--sGS requires fewer iterations
than FRB--Jacobi on all five instances, it has longer full-call times
in every case. These results suggest that SHARP can more effectively
translate the problem structure exploited by the sGS preconditioner
into computational gains on the quadratic recourse games.       

\subsection{Multicommodity-Flow LP Comparison}
\label{sec:lp-solver-comparison}

We consider the multicommodity-flow LP~\eqref{eq:flow-lp}
introduced in Section~\ref{subsec:network-schur}.
The underlying network has nodes $0,\ldots,N-1$ and $E=2N-3$ directed edges, consisting of chain edges $(i,i+1)$ for $i=0,\ldots,N-2$ and chord edges $(i,i+2)$ for $i=0,\ldots,N-3$. We use the full node-edge incidence matrix $D\in\R^{N\times E}$, as in Section~\ref{subsec:network-schur}.

For each commodity $i=1,\ldots,K$, we sample a source node uniformly from $\{0,\ldots,N-3\}$ and a sink node uniformly from the nodes with larger indices. The demand magnitude $\theta_i$ is sampled uniformly from $[0.5,1.5)$. Let $\widehat f_i$ denote the flow that routes $\theta_i$ units
along the chain path from its source to its sink, and set $d_i=D\widehat f_i$. The edge capacities are generated as
\[
a_e^{(0)}
=\sum_{i=1}^K(\widehat f_i)_e+0.05(0.5+V_e),
\qquad
a_e=
\begin{cases}
a_e^{(0)}, & e\text{ is a chain edge},\\
(K/1024)a_e^{(0)}, & e\text{ is a chord edge},
\end{cases}
\]
where the $V_e$ are independent uniform random variables on $[0,1)$. This construction ensures the feasibility of the resulting LP.

All commodities share the same edge-cost vector
$\widetilde c\in\R^E$, with $\widetilde c_e=1$ for chain edges
and $\widetilde c_e=1.875-0.125W_e$ for chord edges,
where the $W_e$ are independent uniform random variables on
$[0,1)$.
We consider two sequences of instances, $K=2N$ and $N=2K$,
each involving $(K+1)E$ primal variables.

For SHARP, we use the degenerate block preconditioner $P_\beta$ developed in Section~\ref{subsec:network-schur}. The corresponding weighted Lipschitz bound is $L_{P_\beta}=1$ for every fixed $\beta>0$, and we choose the endpoint step-size $\gamma=1/2$. The generalized backward step is computed using the Schur-complement formula~\eqref{eq:flow-arrow-solve}, with the zero-mean Laplacian systems solved through a banded factorization of a positive-definite grounded principal minor of $DD^\top$. This factorization is independent of $\beta$ and can be reused across commodities, iterations, and restarts.

We initialize $\beta=1$ and keep it fixed within each restart cycle. At a restart, let $(\Delta x,\Delta y)$ denote the displacement from the current anchor to the selected restart point. Motivated by the $P$-distance term in~\eqref{eq:darg-rate}, we use this displacement as an estimate of the displacement
from the anchor to a solution. To balance the primal and dual contributions to $\|(\Delta x,\Delta y)\|_{P_\beta}^2$, we minimize $\beta\|\Delta x\|^2+\beta^{-1}\|A^\top\Delta y\|^2$
over $\beta\in[10^{-8},10^8]$ and obtain
\[
\beta_{\mathrm{new}}
=
\min\left\{10^8,\,
\max\left\{10^{-8},
\frac{\|A^\top\Delta y\|}{\|\Delta x\|}
\right\}\right\}.
\]
If either $\|\Delta x\|$ or $\|A^\top\Delta y\|$ is
numerically zero, or if the resulting ratio is not finite, we retain the current value of $\beta$.

We compare SHARP with the official HPR-LP-C implementation~\cite{chen2026hprlp} and the MOSEK 11.2.4 interior-point solver~\cite{mosek2026}. SHARP starts from zero, while HPR-LP-C and MOSEK use their default initializations. 

For each returned primal--dual point $(x,y,z)$, the primal and dual objective values for the original LP are given by
\[
p_{\mathrm{obj}}=c^\top x,
\qquad
d_{\mathrm{obj}}=b^\top y+u^\top z_-,
\]
where $z_-=\min\{z,0\}$ componentwise. The relative KKT residuals including primal feasibility, dual stationarity, and the primal--dual gap are defined by
\begin{equation*}\label{eq:flow-common-residual}
r_p=
\frac{\max\{\|Ax-b\|,\|x-\Pi_{[0,u]}(x)\|\}}{1+\|b\|},
\quad
r_d=
\frac{\|c-A^\top y-z\|}{1+\|c\|},
\quad
r_g=
\frac{|p_{\mathrm{obj}}-d_{\mathrm{obj}}|}
{1+|p_{\mathrm{obj}}|+|d_{\mathrm{obj}}|}.
\end{equation*}
For SHARP, we recover $z=c-A^\top y$, which enforces dual
stationarity, and terminate when
$\max\{r_p,r_d,r_g\}\le10^{-6}$.
The comparison solvers retain their native stopping rules. We report full-call times after warm-up, including model
construction, solver setup, iterations, and solution recovery, but excluding instance generation and warm-up.

\begin{table}[!htbp]
\centering
\small
\setlength{\tabcolsep}{2.8pt}
\renewcommand{\arraystretch}{1.16}
\caption{Multicommodity-flow LP comparison: full-call
times (s) after warm-up and iteration counts.
$\mathrm{T}$ denotes the $3600$-second native time limit,
and $\mathrm{M}$ denotes memory exhaustion.}
\label{tab:current-flow}
\vspace{4pt}
\begin{tabular*}{\textwidth}{@{\extracolsep{\fill}}rrrrrrr@{}}
\toprule
Problem size
& \multicolumn{2}{c}{SHARP}
& \multicolumn{2}{c}{HPR-LP-C}
& \multicolumn{2}{c}{MOSEK}\\
\cmidrule(lr){1-1}
\cmidrule(lr){2-3}
\cmidrule(lr){4-5}
\cmidrule(lr){6-7}
$(N,K)$
& Time (s) & Iterations
& Time (s) & Iterations
& Time (s) & Iterations\\
\midrule
$(1\,024,2\,048)$
& \textbf{0.96} & 1\,950
& 66.05 & 189\,300
& 11.97 & 0\\
$(2\,048,4\,096)$
& \textbf{5.22} & 2\,250
& 438.29 & 693\,150
& 66.04 & 8\\
$(4\,096,8\,192)$
& \textbf{24.44} & 2\,400
& $\mathrm{T}$ & ---
& 468.60 & 10\\
$(8\,192,16\,384)$
& \textbf{111.40} & 2\,700
& $\mathrm{M}$ & ---
& $\mathrm{M}$ & ---\\
\midrule
$(2\,048,1\,024)$
& \textbf{0.95} & 1\,800
& 82.44 & 309\,150
& 12.76 & 0\\
$(4\,096,2\,048)$
& \textbf{4.18} & 1\,650
& 530.70 & 831\,150
& 49.87 & 0\\
$(8\,192,4\,096)$
& \textbf{21.39} & 1\,950
& $\mathrm{T}$ & ---
& 276.93 & 9\\
$(16\,384,8\,192)$
& \textbf{94.37} & 2\,100
& $\mathrm{M}$ & ---
& $\mathrm{M}$ & ---\\
\bottomrule
\end{tabular*}
\end{table}

Under the reported configurations, Table~\ref{tab:current-flow} shows that
SHARP solves all eight instances, compared with four for HPR-LP-C and six
for MOSEK. MOSEK solves three of these instances in presolve, with zero
interior-point iterations. On the instances solved by each comparator,
SHARP achieves speedups of approximately $68.9$--$127.0$ over HPR-LP-C
and $11.9$--$19.2$ over MOSEK. Moreover, only SHARP solves the two
largest instances, each containing approximately $2.68\times10^8$
variables, for which HPR-LP-C and MOSEK reach the GPU and host memory
limits, respectively.

These results highlight the benefits of exploiting the shared network
structure through the degenerate block preconditioner. The simplified
backward step, reusable factorization, and matrix-free operator
evaluations reduce computational and storage costs, helping SHARP
solve the largest tested instances within the available memory.

\section{Conclusions}\label{sec:conclusion}

In this paper, we developed SHARP, an accelerated reflected
splitting method for monotone inclusions that accommodates
problem-tailored, possibly degenerate preconditioners.
A reduced-space analysis yields an explicit $\mathcal O(1/k)$
last-iterate natural-residual bound and convergence in the
preconditioner-induced seminorm, with strong convergence in
the original space under an additional continuity condition.
For affine forward operators in finite dimensions, an equivalent
product-space Halpern representation provides convergence
guarantees under alternative assumptions, which relax the
forward range requirement and permit equality in the step-size
condition. Building on these results, we constructed symmetric
Gauss--Seidel and block Schur preconditioners to simplify the
generalized backward step. Numerical experiments demonstrated the computational benefits of combining SHARP with problem-tailored preconditioners on the tested problem families.
These results motivate the development of splitting algorithms and preconditioners tailored to other classes of structured
monotone inclusions, which we leave for future work.

{\sloppy
\bibliographystyle{plain}
\bibliography{ref}

@inproceedings{golowich2020tight,
  author    = {Golowich, Noah and Pattathil, Sarath and Daskalakis, Constantinos},
  title     = {Tight Last-Iterate Convergence Rates for No-Regret Learning in Multi-Player Games},
  booktitle = {Advances in Neural Information Processing Systems},
  volume    = {33},
  pages     = {20766--20778},
  year      = {2020}
}

@book{rockafellar1998variational,
  author    = {Rockafellar, R. Tyrrell and Wets, Roger J.-B.},
  title     = {Variational Analysis},
  publisher = {Springer},
  year      = {1998}
}

@inproceedings{applegate2021practical,
  author    = {David Applegate and Mateo D{\'i}az and Oliver Hinder and Haihao Lu and Miles Lubin and Brendan O'Donoghue and Warren Schudy},
  title     = {Practical Large-Scale Linear Programming using Primal-Dual Hybrid Gradient},
  booktitle = {Advances in Neural Information Processing Systems},
  volume    = {34},
  pages     = {20243--20257},
  year      = {2021},
  note      = {\href{https://proceedings.neurips.cc/paper/2021/hash/a8fbbd3b11424ce032ba813493d95ad7-Abstract.html}{NeurIPS 2021}}
}

@book{bauschke2017convex,
  author    = {Heinz H. Bauschke and Patrick L. Combettes},
  title     = {Convex Analysis and Monotone Operator Theory in Hilbert Spaces},
  edition   = {2},
  publisher = {Springer},
  address   = {Cham},
  year      = {2017},
  note      = {\url{https://doi.org/10.1007/978-3-319-48311-5}}
}

@article{chen2019twostagesge,
  author = {Chen, Xiaojun and Shapiro, Alexander and Sun, Hailin},
  title = {Convergence Analysis of Sample Average Approximation of Two-stage Stochastic Generalized Equations},
  journal = {SIAM Journal on Optimization},
  volume = {29},
  number = {1},
  pages = {135--161},
  year = {2019},
  doi = {10.1137/17M1162822},
  url = {https://doi.org/10.1137/17M1162822}
}

@book{facchinei2003finite,
  author    = {Francisco Facchinei and Jong-Shi Pang},
  title     = {Finite-Dimensional Variational Inequalities and
               Complementarity Problems},
  volume    = {I},
  publisher = {Springer},
  address   = {New York},
  year      = {2003},
  note      = {\url{https://doi.org/10.1007/b97543}}
}

@article{bot2025fast,
  author  = {Radu Ioan Bo{\c{t}} and Dang-Khoa Nguyen and Chunxiang Zong},
  title   = {Fast Reflected Forward--Backward Algorithm: Achieving Fast Convergence Rates for Convex Optimization with Linear Cone Constraints},
  journal = {Journal of Scientific Computing},
  volume  = {105},
  number  = {3},
  pages   = {73},
  year    = {2025},
  note    = {\url{https://doi.org/10.1007/s10915-025-03103-9}}
}

@article{bredies2022degenerate,
  author  = {Kristian Bredies and Enis Chenchene and Dirk A. Lorenz and Emanuele Naldi},
  title   = {Degenerate Preconditioned Proximal Point Algorithms},
  journal = {SIAM Journal on Optimization},
  volume  = {32},
  number  = {3},
  pages   = {2376--2401},
  year    = {2022},
  note    = {\url{https://doi.org/10.1137/21M1448112}}
}

@article{burke1999variable,
  author  = {James V. Burke and Maijian Qian},
  title   = {A Variable Metric Proximal Point Algorithm for Monotone Operators},
  journal = {SIAM Journal on Control and Optimization},
  volume  = {37},
  number  = {2},
  pages   = {353--375},
  year    = {1999},
  note    = {\url{https://doi.org/10.1137/S0363012992235547}}
}

@inproceedings{cai2023accelerated,
  author    = {Yang Cai and Weiqiang Zheng},
  title     = {Accelerated Single-Call Methods for Constrained Min-Max Optimization},
  booktitle = {The Eleventh International Conference on Learning Representations},
  year      = {2023},
  note      = {\href{https://openreview.net/forum?id=HRwN7IQLUKA}{OpenReview:HRwN7IQLUKA}}
}

@article{cai2016zerosum,
  author  = {Yang Cai and Ozan Candogan and Constantinos Daskalakis and Christos Papadimitriou},
  title   = {Zero-Sum Polymatrix Games: A Generalization of Minmax},
  journal = {Mathematics of Operations Research},
  volume  = {41},
  number  = {2},
  pages   = {648--655},
  year    = {2016},
  note    = {\url{https://doi.org/10.1287/moor.2015.0745}}
}

@article{halpern1967fixed,
  author  = {Benjamin Halpern},
  title   = {Fixed Points of Nonexpanding Maps},
  journal = {Bulletin of the American Mathematical Society},
  volume  = {73},
  number  = {6},
  pages   = {957--961},
  year    = {1967},
  note    = {\url{https://doi.org/10.1090/S0002-9904-1967-11864-0}}
}

@article{kim2021accelerated,
  author  = {Donghwan Kim},
  title   = {Accelerated Proximal Point Method for Maximally Monotone Operators},
  journal = {Mathematical Programming},
  volume  = {190},
  pages   = {57--87},
  year    = {2021},
  note    = {\url{https://doi.org/10.1007/s10107-021-01643-0}}
}

@article{korpelevich1976extragradient,
  author  = {Galina M. Korpelevich},
  title   = {The Extragradient Method for Finding Saddle Points and Other Problems},
  journal = {Ekonomika i Matematicheskie Metody},
  volume  = {12},
  number  = {4},
  pages   = {747--756},
  year    = {1976},
  note    = {\url{https://cs.uwaterloo.ca/~y328yu/classics/extragrad.pdf}}
}

@article{rockafellar1976monotone,
  author  = {R. Tyrrell Rockafellar},
  title   = {Monotone Operators and the Proximal Point Algorithm},
  journal = {SIAM Journal on Control and Optimization},
  volume  = {14},
  number  = {5},
  pages   = {877--898},
  year    = {1976},
  note    = {\url{https://doi.org/10.1137/0314056}}
}

@article{sabach2017first,
  author  = {Shoham Sabach and Shimrit Shtern},
  title   = {A First Order Method for Solving Convex Bilevel Optimization Problems},
  journal = {SIAM Journal on Optimization},
  volume  = {27},
  number  = {2},
  pages   = {640--660},
  year    = {2017},
  note    = {\url{https://doi.org/10.1137/16M105592X}}
}

@article{lieder2021convergence,
  author  = {Felix Lieder},
  title   = {On the Convergence Rate of the {Halpern}-Iteration},
  journal = {Optimization Letters},
  volume  = {15},
  pages   = {405--418},
  year    = {2021},
  note    = {\url{https://doi.org/10.1007/s11590-020-01617-9}}
}

@article{li2019sgs,
  author  = {Xudong Li and Defeng Sun and Kim-Chuan Toh},
  title   = {A Block Symmetric {Gauss--Seidel} Decomposition Theorem for Convex Composite Quadratic Programming and Its Applications},
  journal = {Mathematical Programming},
  volume  = {175},
  number  = {1--2},
  pages   = {395--418},
  year    = {2019},
  note    = {\url{https://doi.org/10.1007/s10107-018-1247-7}}
}

@article{malitsky2020forward,
  author  = {Yura Malitsky and Matthew K. Tam},
  title   = {A Forward--Backward Splitting Method for Monotone Inclusions Without Cocoercivity},
  journal = {SIAM Journal on Optimization},
  volume  = {30},
  number  = {2},
  pages   = {1451--1472},
  year    = {2020},
  note    = {\url{https://doi.org/10.1137/18M1207260}}
}

@inproceedings{cai2024accelerated,
  author    = {Yang Cai and Argyris Oikonomou and Weiqiang Zheng},
  title     = {Accelerated Algorithms for Constrained Nonconvex--Nonconcave Min--Max Optimization and Comonotone Inclusion},
  booktitle = {Proceedings of the 41st International Conference on Machine Learning},
  series    = {Proceedings of Machine Learning Research},
  volume    = {235},
  pages     = {5312--5347},
  year      = {2024},
  note      = {\href{https://proceedings.mlr.press/v235/cai24f.html}{PMLR 235} }
}

@inproceedings{lee2021fast,
  author    = {Sucheol Lee and Donghwan Kim},
  title     = {Fast Extra Gradient Methods for Smooth Structured Nonconvex--Nonconcave Minimax Problems},
  booktitle = {Advances in Neural Information Processing Systems},
  volume    = {34},
  pages     = {22588--22600},
  year      = {2021},
  note      = {\href{https://proceedings.neurips.cc/paper/2021/hash/be767243ca8f574c740fb4c26cc6dceb-Abstract.html}{NeurIPS 2021}}
}

@article{tseng2000modified,
  author  = {Paul Tseng},
  title   = {A Modified Forward--Backward Splitting Method for Maximal Monotone Mappings},
  journal = {SIAM Journal on Control and Optimization},
  volume  = {38},
  number  = {2},
  pages   = {431--446},
  year    = {2000},
  note    = {\url{https://doi.org/10.1137/S0363012998338806}}
}

@misc{trandinh2025accelerated,
  author        = {Quoc Tran-Dinh and Nghia Nguyen-Trung},
  title         = {Accelerated Extragradient-Type Methods: Generalization and Unified Convergence Analysis},
  year          = {2025},
  eprint        = {2501.04585},
  archiveprefix = {arXiv},
  primaryclass  = {math.OC},
  note          = {\href{https://arxiv.org/abs/2501.04585v2}{arXiv:2501.04585v2}}
}

@article{sun2025accelerating,
  author  = {Defeng Sun and Yancheng Yuan and Guojun Zhang and Xinyuan Zhao},
  title   = {Accelerating Preconditioned {ADMM} via Degenerate Proximal Point Mappings},
  journal = {SIAM Journal on Optimization},
  volume  = {35},
  number  = {2},
  pages   = {1165--1193},
  year    = {2025},
  note    = {\url{https://doi.org/10.1137/24M1650053}}
}

@article{chen2026hprlp,
  author  = {Kaihuang Chen and Defeng Sun and Yancheng Yuan and Guojun Zhang and Xinyuan Zhao},
  title   = {{HPR-LP}: An Implementation of an {HPR} Method for Solving Linear Programming},
  journal = {Mathematical Programming Computation},
  volume  = {18},
  pages   = {183--210},
  year    = {2026},
  note    = {\url{https://doi.org/10.1007/s12532-025-00292-0}}
}

@misc{mittelmann2026lpfeas,
  author = {Hans D. Mittelmann},
  title  = {{LPfeas} Benchmark (find {PD} feasible point) and Addendum},
  year   = {2026},
  note   = {\url{https://plato.asu.edu/ftp/lpfeas.html}. Benchmark dated September 16, 2026}
}

@manual{mosek2026,
  author = {{MOSEK ApS}},
  title  = {{MOSEK Optimizer API for Python 11.2.4}},
  year   = {2026},
  note   = {\href{https://docs.mosek.com/latest/pythonapi/index.html}{Version 11.2.4}}
}

@misc{sun2026sgs,
  author       = {Defeng Sun and Hailing Wang and Wei Zhao},
  title        = {A Symmetric {Gauss--Seidel} Based Alternating Proximal {ALM} for Generalized {Nash} Equilibrium Problems in {Banach} Spaces},
  year         = {2026},
  eprint       = {2606.13198},
  archiveprefix = {arXiv},
  primaryclass = {math.OC},
  note         = {\href{https://arxiv.org/abs/2606.13198}{arXiv:2606.13198}}
}

@article{yoon2025flock,
  author  = {TaeHo Yoon and Ernest K. Ryu},
  title   = {Accelerated Minimax Algorithms Flock Together},
  journal = {SIAM Journal on Optimization},
  volume  = {35},
  number  = {1},
  pages   = {180--209},
  year    = {2025},
  note    = {\url{https://doi.org/10.1137/22M1504597}}
}

@article{wittmann1992approximation,
 author={Rainer Wittmann}, title={Approximation of Fixed Points of Nonexpansive Mappings},
 journal={Archiv der Mathematik}, volume={58}, pages={486--491}, year={1992},
 note={\url{https://doi.org/10.1007/BF01190119}}
}

@article{cevher2021reflected,
 author={Volkan Cevher and B{\`a}ng C{\^o}ng V{\~u}},
 title={A Reflected Forward--Backward Splitting Method for Monotone Inclusions Involving {Lipschitzian} Operators},
 journal={Set-Valued and Variational Analysis}, volume={29}, pages={163--174}, year={2021},
 note={\url{https://doi.org/10.1007/s11228-020-00542-4}}
}

@article{borgens2018distributed,
  title={A distributed regularized Jacobi-type ADMM-method for generalized Nash equilibrium problems in Hilbert spaces},
  author={B{\"o}rgens, Eike and Kanzow, Christian},
  journal={Numerical Functional Analysis and Optimization},
  volume={39},
  number={12},
  pages={1316--1349},
  year={2018},
  publisher={Taylor \& Francis}
}
}

\end{document}